\documentclass[a4paper,11pt]{article}

\usepackage[utf8]{inputenc}   
\usepackage[T1]{fontenc}      
\usepackage{geometry}         
\usepackage[francais,english]{babel}  

\newtheorem{theorem}{Theorem}
\newtheorem{proposition}[theorem]{Proposition}
\newtheorem{lemma}[theorem]{Lemma}
\newtheorem{corollary}[theorem]{Corollary}
\newtheorem{definition}[theorem]{Definition}

\newtheorem{example}[theorem]{Example}
\newtheorem{notation}[theorem]{Notation}

\newtheorem{remark}[theorem]{Remark}

\newcommand{\qed}{\nobreak \ifvmode \relax \else
      \ifdim\lastskip<1.5em \hskip-\lastskip
      \hskip1.5em plus0em minus0.5em \fi \nobreak
      \vrule height0.75em width0.5em depth0.25em\fi}

\makeindex 

\usepackage{pb-diagram}  
\usepackage{makeidx}     
\usepackage{graphicx}    
\usepackage{multicol}    
\usepackage{cite}        
\usepackage{amsfonts}
\usepackage{amssymb}
\usepackage{amscd}
\usepackage{epsfig}

\def\C{\mathbf{C}}
\def\F{\mathbf{F}}
\def\T{\mathbf{T}}

\def\Z{\mathbf{Z}} 
\def\R{\mathbf{R}} 
\def\N{\mathbf{N}}

\def\G{\mathbf{G}}    

\def\a{\alpha}

\def\l{\lambda}

\def\rmi{\uppercase\expandafter{\romannumeral1}}
\def\rmii{\uppercase\expandafter{\romannumeral2}}
\def\rmiii{\uppercase\expandafter{\romannumeral3}}
\def\rmv{\uppercase\expandafter{\romannumeral5}}
\def\rmvi{\uppercase\expandafter{\romannumeral6}}
\def\rmvii{\uppercase\expandafter{\romannumeral7}}
\def\rmviii{\uppercase\expandafter{\romannumeral8}}

\def\RR{\mathcal{R}}              

\def\FF{\mathcal{F}}          
\def\TT{\mathcal{T}}

\def\TP{{\mathcal{T}}^2}

\def\I{\mathcal{I}}
\def\J{\mathcal{J}}              

\def\X{\mathcal{X}}              

\def\g{\mathfrak{g}}          
\def\z{\mathfrak{z}}
\def\gg{{\mathfrak{g}}^2}      
\def\h{\mathfrak{h}}          
\def\gp{\g_+}                 
\def\gm{\g_-}                 
\def\ggm{{\g}^2_-}             
\def\ggp{{\g}^2_+}             
\def\Liesl#1{\mathfrak{sl}_{#1}}     
\def\gl#1{\mathfrak{gl}_{#1}}        

\def\can#1{\left\langle#1\right\rangle}

\def\PB{\left\{\cdot\,,\cdot\right\}}

\def\Pb#1{\left\{\cdot\,,#1\right\}}
\def\pb#1{\left\{#1\right\}}

\def\lb#1{\[#1\]}
\def\LB{[\cdot\,,\cdot]}

\def\inn#1#2{\left\langle#1\,\vert\,#2\right\rangle}

\def\INN{\langle\cdot\,\vert\cdot\,\rangle}

\def\({\left(}
\def\){\right)}
\def\[{\left[}
\def\]{\right]}

\def\ad{\mathop{\rm ad}\nolimits}
\def\Ad{\mathop{\rm Ad}\nolimits}

\def\cycl{\mathop{\rm cycl}\nolimits}

\def\dim{\mathop{\rm dim}\nolimits}
\def\diff{\mathsf{d}}

\newcommand{\Rk}{\mathop{\rm Rk}\nolimits}

\def\card{\mathop{\rm card}\nolimits}

\def\Trace{\mathop{\rm Trace}\nolimits}

\newenvironment{eqn*}[1][1.5]
  {$$\renewcommand{\arraystretch}{#1}
      \begin{array}{rcl}}
      {\end{array}$$}

\def\ds{\displaystyle}
\def\comment#1{}  

\def\p{\partial}
\def\pp#1#2{\frac{\p #1}{\p #2}}

\catcode`\á=13 \defá{\'a}
\catcode`\à=13 \defà{\`a}
\catcode`\â=13 \defâ{\^a}
\catcode`\ä=13 \defä{\"a}
\catcode`\é=13 \defé{\'e}
\catcode`\è=13 \defè{\`e}
\catcode`\ê=13 \defê{\^e}
\catcode`\ë=13 \defë{\"e}
\catcode`\í=13 \defí{\'\i}
\catcode`\ì=13 \defì{\`\i}
\catcode`\î=13 \defî{\^\i}
\catcode`\ï=13 \defï{\"\i}
\catcode`\ó=13 \defó{\'o}
\catcode`\ò=13 \defò{\`o}
\catcode`\ô=13 \defô{\^o}
\catcode`\ö=13 \defö{\"o}
\catcode`\ú=13 \defú{\'u}
\catcode`\ù=13 \defù{\`u}
\catcode`\û=13 \defû{\^u}
\catcode`\ü=13 \defü{\"u}

\catcode`\ç=13\defç{\c c}

\catcode`\Á=13 \defÁ{\'A}
\catcode`\À=13 \defÀ{\`A}
\catcode`\Â=13 \defÂ{\^A}
\catcode`\Ä=13 \defÄ{\"A}
\catcode`\É=13 \defÉ{\'E}
\catcode`\È=13 \defÈ{\`E}
\catcode`\Ê=13 \defÊ{\^E}
\catcode`\Ë=13 \defË{\"E}
\catcode`\Í=13 \defÍ{\'I}
\catcode`\Ì=13 \defÌ{\`I}
\catcode`\Î=13 \defÎ{\^I}
\catcode`\Ï=13 \defÏ{\"I}
\catcode`\Ó=13 \defÓ{\'O}
\catcode`\Ò=13 \defÒ{\`O}
\catcode`\Ô=13 \defÔ{\^O}
\catcode`\Ö=13 \defÖ{\"O}
\catcode`\Ú=13 \defÚ{\'U}
\catcode`\Ù=13 \defÙ{\`U}
\catcode`\Û=13 \defÛ{\^U}
\catcode`\Ü=13 \defÜ{\"U}

\catcode`\Ç=13\defÇ{\C C}

\begin{document}
\author{Khaoula Ben Abdeljelil}
\title{The integrability of the $2$-Toda lattice on  a simple Lie algebra}
 \date{}                     
\maketitle
\begin{abstract}
 We define the $2$-Toda lattice on every simple Lie algebra $\g$,  and we show its Liouville integrability. 
We show that this  lattice is given by a pair of Hamiltonian vector fields, associated with a Poisson bracket 
which results from an $\RR$-matrix of the 
underlying Lie algebra. We construct a big family of constants
 of motion which  we use to prove the Liouville integrability of the system.
 We achieve  the proof of their integrability by using several results on simple Lie algebras,
 $\RR$-matrices, invariant functions and root systems.  
\end{abstract}

\tableofcontents

\section{Introduction}
  The  \emph{$2$-Toda lattice associated with $\Liesl{n}(\C)$} is the
  pair of differential  equations given
 by the Lax equations
 \begin{equation}\label{k1}
\renewcommand{\arraystretch}{2.5}
\left\{
\begin{array}{rcl}
\ds\pp{(L,M)}{t}&=&[{(L_u, L_u)},(L,M)],\\
\ds\pp{(L,M)}{s}&=&[{(M_{\ell},M_{\ell})},(L,M)],
\end{array}
\right.
\end{equation}%
where  $(L,M)$ are traceless matrices of the form 
\begin{equation}\label{es}
(L,M)=
\left(
\begin{array}{c}
\left(
\begin{array}{cccc}
a_{11}&1     &         &0     \\
a_{21}&a_{22}&\ddots   &      \\
\vdots&      &\ddots   &1     \\
a_{n1}&\cdots&a_{n,n-1}&a_{nn}
\end{array}
\right)
,
\left(
\begin{array}{cccc}
b_{11}&\cdots&\cdots   &b_{1n}\\
b_{21}&\ddots&         &\vdots\\
      &\ddots&\ddots   &\vdots\\
0     &      &b_{n,n-1}&b_{nn}
\end{array}
\right)
\end{array}
\right),
\end{equation}%
and where $L_u$ is the upper triangular part of $L$ and $M_{\ell}$ is the
strictly lower triangular part of $M$.

The  $2$-Toda lattice  was first introduced in the context of infinite-dimensional matrices
\cite{adlerv}, \cite{Carlet}, \cite{ueta}. This system is Hamiltonian with respect to a Poisson structure,
 associated with an $\RR$-matrix with $\frac 12\Trace L^2$ and $\frac 12 \Trace M^2$ as Hamiltonians. 
These Hamiltonians admit the $\Ad$-invariant functions on $\gl{}((\infty))\times\gl{}((\infty))$ 
as Poisson commuting constants of motion. In the finite-dimensional setting,
i.e., on $\Liesl{n}(\C)$, this system is still Hamiltonian with respect to a Poisson structure associated with an  $\RR$-matrix similar to that of  the infinite-dimensional setting. The family of $\Ad$-invariant functions
are again in involution with respect to the Poisson structure associated with the $R$-matrix but their number of independent
$\Ad$-invariant functions, which is $2n-2$, is much too small compared with the dimension  of the phase space
  of the $2$-Toda lattice, which is $n^2+2n-3$.

 The main purpose of the present article is to introduce a family
of functions which contains the $\Ad$-invariant functions  and is large enough to 
give the Liouville integrability of the  $2$-Toda lattice. This family is
$\FF:=(F_{j,i}, 1\leqslant i\leqslant n-1 \textrm{ and } 0\leqslant j \leqslant i+1)$, where $F_{j,i}$ 
is the function on $\Liesl{n}(\C)$ constructed from the  the relation
$$\Trace(\l L-M)^{i+1}=\sum_{j=0}^{i+1}(-1)^{m_i+1-j}\l^j F_{j,i}(L,M).$$
The Toda lattice  can be defined for  every simple Lie algebra \cite{kostant}. In the same spirit  
we introduce the 2-Toda system for every simple Lie algebra and we show its Liouville
 integrability. We recall that, according to  \cite[Definition 4.13]{Liv}, in the general case a system $(M,\PB,\FF)$ is Liouville 
integrable if $(M,\PB)$ is a Poisson manifold of  rank
$2r$ and $\FF=(F_1,\dots,F_s)$ is involutive and independent, with $s=\dim M-r$.  

\bigskip

The main result of this paper is given in Section $3$:  we show the Liouville integrability of the $2$-Toda lattice,
 not only on $\Liesl{n}(\C)$, but also on an arbitrary simple Lie algebra $\g$. We define this lattice
 in this general context.
We choose a simple Lie algebra $\g$; then we denote by $\h$ a Cartan subalgebra, by $
\a_1,\dots,\a_{\ell}$ the simple roots of $\g$ with respect to $\h$, by $e_1,\dots,e_{\ell}$ the corresponding eigenvectors,
and by $\g_i$ the subspace spanned by the eigenvectors associated with roots of length $i \in {\mathbb Z}$.
The \emph{$2$-Toda lattice,
 associated with  $\g$},
 is the  pair of differential equations
on $\TP$ given by  the following Lax pair equations
\begin{equation}\label{equa-motion}
\renewcommand{\arraystretch}{2.5}
\left\{
\begin{array}{ccc}
\renewcommand{\arraystretch}{2.5}
\ds\pp{(L,M)}{t}&=&[(L_+,L_+),(L,M)],\\
\ds\pp{(L,M)}{s}&=&[(M_-,M_-),(L,M)],
\end{array}
\right.
\end{equation}
where  $(L,M)$ belongs to the phase space 
$\TP: = \sum_{i\leqslant 0}\g_i\times\sum_{i\geqslant -1}\g_i+(\sum_{i=1}^{\ell}e_i,0)$, and  
here, $L_+$ (resp.  $M_-$) stands for the projection of $L\in \g$ (resp. $M \in \g$) on 
$\sum_{i\geqslant 0}\g_{i}$ (resp. $\sum_{i<0}\g_i$).

We construct a linear Poisson structure with respect to which the $2$-Toda lattice is 
Hamiltonian. This linear Poisson structure comes from an $\RR$-matrix of the type
studied in Section $2$, namely from the endomorphism of $\gg$ given by $\RR(x,y)=(R(x-y)+y,R(x-y)+x)$,
for all $x,y\in\g$, where
 $R:=P_+-P_- $ is the classical $R$-matrix on a simple Lie algebra, associated with the Lie algebra  splitting
 $\g=\sum_{i\geqslant 0}\g_i\oplus\sum_{i<0} \g_i$. Using the Killing form of $\g$, the product  Lie algebra
 $\gg$ is equipped with a bilinear, symmetric, $\Ad$-invariant, non-degenerate form, which allows us to identify  $\gg$ with its dual and hence to equip $\gg$
with a linear Poisson structure that we denote by $\PB_{\RR}$.
We then establish that the phase space $\TP$ of the $2$-Toda lattice  is a Poisson submanifold of
$(\gg,{\PB}_{\RR})$ (see Propositions \ref{K6}) and the Hamiltonian vector fields of functions 
 $H(x,y):=\frac{1}{2}\inn{x}{x}$ and  ${\tilde{ H}}(x,y):=\frac 12\inn{y}{y}$ are  tangent to $\TP$ and
describe on $\TP$  the equations of  motion~(\ref{equa-motion})  for  the  $2$-Toda lattice (see Proposition \ref{hamil}).

We  construct the integrable system.
  Let $(P_1,\dots,P_{\ell})$ be a generating family of the algebra of $\Ad$-invariant functions,
 chosen to be of respective degrees $m_1+1,\dots,m_{\ell}+1$. For all $ 1\leqslant i\leqslant \ell$ and 
$0\leqslant j\leqslant m_i+1$,
we define $F_{j,i} \in \FF(\gg)$ to be the coefficient in $\l^j$ of the polynomial  $(x,y) \mapsto P_i(\l x-y)$,
 and we consider the following family of functions on $\gg$:
 $$ \FF=(F_{j,i}, 1\leqslant i\leqslant \ell \textrm{ and }0\leqslant j\leqslant m_i+1 ).$$
Our main result is the following theorem:

\bigskip
\noindent
{\bf Theorem}
{\it The  triplet $(\TP,\FF_{\arrowvert\TP},{\PB}_{\RR})$ is  an integrable system.}
\bigskip
\\
To prove this result, we proceed as follows:
\begin{itemize}
\item We show in Proposition \ref{kh6}  that $\FF$ is involutive for ${\PB}_{\RR}$: this result is a particular case of the second point 
of Theorem \ref{tr1} that we prove in Section $2$.
\item We show in Proposition \ref{indtt}  the independence of $\FF$ restricted  to $\TP$ by  using a theorem of Ra\"is \cite{Rais}.
\item We compute the rank of the restriction of the  Poisson $\RR$-bracket to $\TP$.
 This is done  by establishing a Poisson 
isomorphism between $(\TP,{\PB}_{\RR})$ and the product Poisson manifold $(\g^*,\PB)\times(\g_0\oplus\g_1,{\PB}_R)$, 
where $\PB$ is the Lie-Poisson bracket on $\g^*$, and  ${\PB}_R$ is the above $R$-bracket 
(restricted to $\g_0\oplus\g_1$; see Proposition~\ref{isomoTP}).  
\item We check that $\card\FF_{\arrowvert\TP}= \dim\TP-\frac{1}{2}\Rk(\TP,{\PB}_{\RR})$.
\end{itemize}

We also define
 the $2$-Toda from on $\gl{n}(\C)$  and we show  its 
 Liouville integrability  with respect  to the linear and quadratic Poisson $\RR$-bracket. This result is mainly based on the 
fact that the Hamiltonian vector fields of the linear and quadratic Poisson structure are essentially identical.

 In the last subsection,  by restricting the $2$-Toda lattice on $\g$  to a well chosen affine subspace, 
we find again the usual Toda lattice.
\section{$\RR$-matrices, Poisson structures and  functions in involution
  on the square of a Lie algebra}

In this section we  fix  $(\g,\LB)$  a Lie
algebra over  $\F$, with $\F=\R$ or $\F=\C$.
The vector space $\gg:=\g\times\g$, endowed with  the Lie bracket
\begin{equation}\label{equ-crochet -lie}
[(x,y),(z,s)]:=([x,z],[y,s]),\qquad \qquad \forall x,y,z,s\in\g,
\end{equation}%
 is  a Lie algebra.


 We construct an $\RR$-matrix of  $\gg$ with the help  of  an endomorphism $R$ of $\g$ satisfying some  conditions.
 With this $\RR$-matrix, when $\g$ is  finite dimensional,  we  construct
a linear Poisson bracket on the dual  ${\gg}^*$ of $\gg$, we  explain
the construction of a large family
 of functions on ${\gg}^*$ which commute for this linear Poisson structure
and we spell out  the expressions of their Hamiltonian vector
fields. We also  give  some
Casimir  functions.
\subsection{Construction of $\RR$-matrices on $\gg$}\label{parg1}

We begin by recalling some properties and  definitions of  $R$-matrices (see \cite[Section 4.4]{Liv}).
Let  $\g$ be a Lie algebra.  A (vector space)  endomorphism $R$  of   $\g$ is called 
an \emph{$R$-matrix} of  $\g$ if  the bilinear map ${\LB}_R:\g\times\g\to\g$, defined for all $x,y\in\g$ by
\begin{equation}\label{T1}
{[x,y]}_R:=\frac{1}{2}(\lbrack Rx,y\rbrack +\lbrack x,Ry\rbrack ),
\end{equation}%
defines a (second) Lie bracket on $\g$, which is then called a Lie  $R$-bracket.
Let
  $B_R:\g\times\g\to\g$ be defined, for all $x,y\in\g$, by
\begin{equation}\label{T2}
B_R(x,y):=\lbrack Rx,Ry\rbrack -R(\lbrack Rx,y\rbrack +\lbrack
x,Ry\rbrack ).
\end{equation}%
 The bracket   ${\LB}_R$ satisfies the Jacobi identity if and only if
\begin{equation}\label{T3}
\lbrack B_R(x,y),z\rbrack +\lbrack B_R(y,z),x\rbrack +\lbrack
B_R(z,x),y\rbrack =0, \qquad \forall  x,y,z\in\g.
\end{equation}%
%
%
A sufficient condition for (\ref{T3}) to be satisfied is  that there
exists a constant $c \in\F$ such that
\begin{equation}\label{M1}
B_R(x,y)=-c^2[x,y],  ~~\qquad\qquad \forall x,y\in\g.
\end{equation}%
We call  (\ref{M1}) \emph{the modified classical   Yang-Baxter equation}
(mCYBE) of $\g$ of  constant $c$.

We construct in the following proposition an $\RR$-matrix of $\gg$ by using an endomorphism of $\g$ which satisfies
 some condition which generalizes the modified classical   Yang-Baxter equation
 of $\gg$.
\begin{proposition}\label{olm1m}
Let $\g$ be  a Lie algebra,  $\z(\g)$ its  center and $c\in\F$ a
constant. Let $R$ be an endomorphism of $\g$ and $\RR$ the endomorphism
of $\gg$ defined, for every $(x,y)\in\gg$,  by
\begin{equation}\label{ml23p}
{\RR}(x,y):=\left(R(x-y)+cy,R(x-y)+cx\right).
\end{equation}%
  The endomorphism $\RR$
is an  $\RR$-matrix of  $\gg$  if and only if,  for every  $x,y\in\g$,
\begin{equation}\label{condition}
B_R(x,y)+c^2[x,y]\in\z(\g).
\end{equation}%
In the particular case where $\g$ is a  complex semi-simple Lie
algebra, this is in turn equivalent to $R$ being  a  solution of  (mCYBE) of  constant
  $c$.
\end{proposition}
\begin{proof}
Using  Definition (\ref{T2}) of $B_{\RR}$ and replacing $\RR$ by its Expression (\ref{ml23p}), we obtain
\begin{equation}\label{Brr}
\begin{array}{ccc}
B_{\RR}((x,y),(z,s))&=&( B_R(x-y,z-s)+c^2([x-y,z-s]-[x,z]),\\
                     &&B_R(x-y,z-s)+c^2([x-y,z-s]-[y,s]) ).
\end{array}
\end{equation}%
 According to Formula (\ref{T3}),    $\RR$ is an
 $\RR$-matrix  of $\gg$ if and only if,  for every
 $(x,y),(x',y')$ and  $(x'',y'')$ in $\gg$,
\begin{equation}\label{r26p}
[B_{\RR}((x,y),(x',y')),(x'',y'')]+\circlearrowleft=(0,0),
\end{equation}%
where
  $\circlearrowleft=\cycl((x,y),(x',y'),(x'',y''))$.
By  Formula (\ref {Brr}),  the left hand side  of  Equation~(\ref{r26p}) is equal to
\begin{equation}\label{kj1m}
\left([B_{R}(x-y,x'-y')+c^2[x-y,x'-y'],x''],
\lbrack B_{R}(x-y,x'-y')+c^2\lbrack x-y,x'-y'\rbrack ,y''\rbrack \right)+\circlearrowleft,
\end{equation}
where we used the Jacobi identities $
  \lbrack\lbrack x,x'\rbrack,x''\rbrack +\circlearrowleft=0 $ and
$\lbrack\lbrack y,y'\rbrack,y''\rbrack +\circlearrowleft=0$. Now,
if  $\RR$
is an  $\RR$-matrix of  $\gg$, then (\ref {kj1m}) is equal to zero  (even
without circular permutation), so  Formula (\ref{condition}) holds. If Formula  (\ref{condition}) holds then
  the first component  of the
 left hand side   of  (\ref {kj1m}) is equal to zero:
\begin{equation}\label{bvg0}
{[B_{R}(x-y,x'-y')+c^2[x-y,x'-y'],x'']}+\cycl((x,y),(x',y'),(x'',y''))=0,
\end{equation}
for every  $x,y, x',y', x'',y''\in\g$;  in particular, for $x'=x''=0$, we obtain
\begin{eqn*}
\lbrack B_R(y',y'')+c^2\lbrack y',y''\rbrack ,x\rbrack =0, \qquad \forall y',y'',x\in\g,
\end{eqn*}%
hence the endomorphism $\RR$
is an  $\RR$-matrix of  $\gg$ holds.
\end{proof}
For   future use, we also state the following corollary, which is a consequence of 
 Formula~(\ref{Brr}) and  already appeared in \cite{SemenovS} (with $c=1$).
\begin{corollary}\label{pror}
Let  $\g$ be  a Lie algebra  and let  $R$  be an endomorphism of $\g$. If
$R$ is a solution of (mCYBE)  of $\g$  of constant $c$, then the endomorphism
$\RR$ of $\gg$, defined, for
all $(x, y)\in\gg$, by
 \begin{eqn*}
{\RR}(x,y):=(R(x-y)+cy,R(x-y)+cx),
\end{eqn*}%
is a solution of (mCYBE) of $\gg$ of constant $c$.
\end{corollary}
\begin{example}\label{splitting-gg}

Let  $\g:=\gp\oplus\gm$ be a   Lie  algebra splitting (i.e., $\g$ is, as a vector space,
 the direct sum of  Lie subalgebras $\gp$ and $\gm$). Let  $P_{\pm}$  be  the  projections of  $\g$
 on
 $\g_{\pm}$ and let $R:=P_+-P_-$. Since $R$ is a solution of (mCYBE) of $\g$ of $c=1$,
 Corollary \ref{pror} implies that   the endomorphism  $\RR$, defined for every  $x,y\in\g$ 
 by $\RR(x,y)=(R(x-y)+y,R(x-y)+x)$
is a solution of (mCYBE) of $\gg$ of $c=1$. Replacing $R$ by its expression $P_+-P_-$, we verify that
\begin{eqnarray}
\RR(x,y)&=&(x_+-x_-+2y_-,y_--y_++2x_+)\nonumber\\
        &=&(x_++y_-,x_++y_-)-(x_--y_-,y_+-x_+)\label{Cm9},
\end{eqnarray}
where  $x=x_++x_-$ and $y=y_++y_-$ are the decompositions of $x$ and $y$ with respect to the splitting $\g=\gp\oplus\gm$.
Since the Lie subalgebras
\begin{eqn*}
\ggp:=\{(x,x) ~~ | ~~ x\in\g \} ~~\textrm { and } ~~ \ggm:=\{(x,y) ~~ | ~~
x\in\gm\textrm { and } y\in\gp\}
\end{eqn*}%
have a trivial intersection, (\ref{Cm9}) implies that  $\gg=\ggp\oplus\ggm$ is a  Lie algebra splitting and 
that  $\RR$ is  the difference of  the projections on  $\ggp$ and $\ggm$.

\end{example}
Note that independently of Proposition \ref{olm1m} or Corollary \ref{pror} the above example  shows that
 from an $R$-matrix of $\g$,   which is the difference of the  projections on two supplementary  Lie subalgebras of $\g$, 
  the endomorphism $\RR$ given by the Formula (\ref{ml23p}) is automatically an $\RR$-matrix of $\gg$,
 which is the difference of projections on two supplementary Lie  subalgebras of $\gg$.
\subsection{Linear Poisson structures on $\g^*\times\g^*$}\label{14mljhsd}

From now and  until the end of the article,  we assume  that $\g$ is a  finite-dimensional  Lie algebra and we  
denote by   $\FF(\g)$
  the algebra of smooth  functions on $\g$ (if $\g$ is a Lie algebra over $\R$) or holomorphic 
functions on $\g$ (if $\g$ is a  Lie algebra over $\C$).
Let $\RR$ be  an $\RR$-matrix of $\gg$.
 The dual $(\gg)^*$ of the Lie algebra $(\gg,{\LB}_{\RR})$ admits a   Lie-Poisson structure, defined for every  $F,G\in\FF((\gg)^*)$ at
 $\varphi\in(\gg)^*$ by
\begin{eqnarray}\label{pdg2m}
{\pb{F,G}}_{\RR}(\varphi):&=&\can{\varphi,{[\diff_{\varphi}F,\diff_{\varphi}G]}_{\RR}}\\
                          &=&\frac 12\can{\varphi,{[\RR\diff_{\varphi}F,\diff_{\varphi}G]+
                          [\diff_{\varphi}F,\RR\diff_{\varphi}G]}},\nonumber 
\end{eqnarray}%
where $\diff_{\varphi}F\in{((\gg)^*)}^*\backsimeq\gg$ is the differential of $F$ at the point $\varphi\in(\gg)^*$.
 We identify
 $(\gg)^*$ with
  $\g^*\times\g^*$ as follows:
\begin{eqn*}
\begin{array}{ccccc}
\Phi&:&\g^*\times\g^*&\to&(\gg)^*\\
     &&(\xi,\eta)&\mapsto&((x,y)\mapsto \can{(\xi,\eta),(x,y)}):=\can{\xi,x}-\can{\eta,y}.
\end{array}
\end{eqn*}%
 By means of  this identification,  the  Poisson structure
 (\ref{pdg2m})  can be transported to the vector space $\g^* \times\g^* $.   This transported   Poisson
 structure is  given, for all functions $F, G\in\FF (\g^* \times\g^*)$ at   $(\xi, \eta)\in\g^* \times\g^* $, by
\begin{equation}\label{rlpb}
{\pb{F,G}}_{\RR}(\xi,\eta ):=\can{(\xi,\eta),{[\diff_{(\xi,\eta)}
    F,\diff_{(\xi,\eta )}G]}_{\RR}},
\end{equation}%
where $\diff_{(\xi,\eta)}F\in (\g^*\times \g^*)^*\simeq \g\times\g$ . We call this Poisson
 structure   \emph{the Poisson $\RR$-bracket on $\g^*\times\g^*$}.
\subsection{Casimirs and  functions  in involution on $\g^*\times\g^*$}
  
Let $\g$ be a finite-dimensional Lie algebra.
We denote by $\G$ a connected  Lie group whose Lie algebra is $\g$. We denote by  ${\FF(\g^*)}^{\G}$ the algebra of $\Ad^*$-invariant functions on $\g^*$.

According to the Adler-Kostant-Symes theorem \cite [Theorem 4.37]{Liv},
the $\Ad^*$-invariant functions on $\g^*\times\g^*$ (i.e., functions in ${\FF(\g^*\times\g^*)}^{\G\times\G}$) are  in involution for the Poisson
structure associated with an $\RR$-matrix. But, when $\g$ is a semi-simple Lie algebra,  this family of functions  is  too  small  to insure  the Liouville
integrability. We construct a larger  family\footnote{Which contains the $\Ad^*$-invariant functions on $\g^*\times\g^*$.}
 of   functions in involution.
\begin{theorem}\label{tr1}
Let~$\lambda\in\F$ be  a
constant and   let  $\psi_{\lambda}$ be  the map	
\begin{equation}\label{C9}
\begin{array}{ccccc}
\psi _{\lambda}&:&\g^*\times\g^*  &\to    &\g ^* \\
               &&(\xi,\eta )&\mapsto &\lambda \xi-\eta.
\end{array}
\end{equation}%
Let  $R$ be an  endomorphism\footnote{$R$ is  not   necessarily an
  $R$-matrix.} of  $\g$ and  $c\in\F$   a constant. Assume that the  endomorphism
$\RR$ of
$\gg$ defined  for every  $(x,y)\in\gg$ by

\begin{equation}\label{formula1}
  \RR(x,y):=\left(R(x-y)+cy,R(x-y)+cx\right)
\end{equation}
is   an $\RR$-matrix of  $\gg$, and  denote by     ${\PB}_{\RR}$  the
Poisson $\RR$-bracket  on $\g^*\times\g^*$.  Then:
\begin{enumerate}
\item[(1)]For every   $F\in\FF(\g^*)^{\G}$, the
  function     $F\circ \psi_1$  is  a Casimir for   ${\PB}_{\RR}$.
\item[(2)] For every   $F,G\in\FF(\g^*)^{\G}$ and
  every   $\lambda,\gamma\in\F$, the
  functions $F\circ \psi_{\lambda}$ and  $G\circ \psi_{\gamma}$ are
 in  involution for
  ${\PB}_{\RR}$. In particular, if  $F$ and  $G$ are  polynomials  of
  degree respectively  $l$ and  $k$, then  the  functions
   $F_0,\dots,F_l,G_0,\dots,G_k$, defined, for every
  $(\xi,\eta)\in\g^*\times\g^*$, by
\begin{eqn*}
  F(\psi_{\l}(\xi,\eta))=\sum_{i=0}^{l} \lambda^i
  F_i(\xi,\eta)\qquad \textrm{ and  }\qquad
G(\psi_{\gamma}(\xi,\eta))=\sum_{j=0}^{k} \gamma^j
  G_j(\xi,\eta),
\end{eqn*}%
are  in involution  for   ${\PB}_{\RR}$.
\item[(3)]  If  $c=1$,  the map  $\psi _1:
  (\g^*\times\g^*,{\PB}_{\RR})\to(\g ^*,\PB)$ is a  Poisson  morphism.
\item[(4)]  For every    $H\in\FF{(\g^*)}^{\G}$,
the   Hamiltonian  vector field $\X_{H\circ\psi_{\lambda}}:={\Pb{H\circ\psi
    _{\lambda}}}_{\RR}$  is given at      $(\xi,\mu)\in\g^*\times\g^*$ by
\begin{equation}\label{champh}
\X_{H\circ\psi _{\lambda}}(\xi,\eta )=\frac{1}{2}(1-\lambda
)\ad^*_{(( R-cI )\diff_{\lambda \xi-\eta} H,( R+cI )\diff_{\lambda \xi-\eta } H)}(\xi,\eta ).
\end{equation}%
In particular,  if $\g:=\gp\oplus \gm$  is a  Lie algebra splitting  and 
    $R:=P_+-P_-$  is the difference of the  projections on $\gp$ and $\gm$, we have
\begin{equation}\label{relchampspli}
\X_{H\circ\phi _{\lambda}}(\xi,\eta )=(1-\lambda )\ad^*_{((-
  \diff_{\lambda \xi-\eta } H)_-,(\diff_{\lambda \xi-\eta } H )_+) }(\xi,\eta ).
\end{equation}
\end{enumerate}
\end{theorem}
To prove this theorem, we need  the following lemmas, the proof of the first of  which is left to the reader.

\begin{lemma}\label{lemma2}
Let  $F\in\FF(\g^*)$,  $(\xi,\eta)\in{{\g}^*\times\g^*}$ and  $\l\in\F$ be  a
constant.
\begin{enumerate}
\item[(1)] The   differential  of  $F\circ\psi_{\l}$ at $
(\xi,\eta)\in\g^*\times\g^*$  is given by
\begin{equation}\label{T15}
\diff_{(\xi,\eta) }(F\circ\psi _{\lambda})=(\lambda
\diff_{\lambda \xi-\eta }{F},\diff_{\lambda \xi-\eta } {F}).
\end{equation}%
%
\item[(2)]  Let  $R$ be  an endomorphism of  $\g$ and  let $\RR$ be   the endomorphism
of  $\gg$ defined  in~(\ref{formula1}). Then
\begin{equation}\label{TR}
\RR\,\diff_{(\xi,\eta) }(F\circ\psi _{\lambda})=(\lambda-1 )
                                          (R\,\diff_{\lambda \xi-\eta }{F},R\,\diff_{\lambda \xi-\eta }{F})
                                          +c(\diff_{\lambda \xi-\eta }{F},\lambda\diff_{\lambda \xi-\eta }{F} ).
\end{equation}
\end{enumerate}
\end{lemma}
\begin{lemma}\label{commute1}
Let    $F,G\in\FF{(\g^*)}^{\G}$.
\begin{enumerate}
\item[(1)]  The differentials of  $F$ and     $G$  commute at every  point of $\g^*$, i.e.,
\begin{equation}\label{btp}
  [\diff_{\xi}F,\diff_{\xi}G]=0,\qquad \qquad \forall\xi\in\g^*.
\end{equation}%
\item[(2)] For every $\xi,\eta \in\g^*$ and every   $a, b \in \F$,
\begin{equation}\label{T18}
\can{a\xi+b\eta,[\diff_{\xi}{F},\diff_{\eta}{G}]}=0.
\end{equation}
\end{enumerate}
\end{lemma}

\begin{proof}
 Using the $\Ad^*$-invariance of  $F$ and  $G$, we  obtain
\begin{equation}\label{fg09}
\ad^*_{\diff_{\xi}F}\xi=0\qquad\textrm{ and }\qquad
\ad^*_{\diff_{\xi}G}\xi=0.
\end{equation}%
Stated differently, $
\diff_{\xi}F\in\g^{\xi}$ and $\diff_{\xi}G\in\g^{\xi}$ where $\g^{\xi}$ is the centralizer of $\xi$. For   $\xi$   regular, according to
  \cite[Proposition  19.7.5]{Yu},
  the  centralizer   $\g^{\xi}$ is  abelian,  so Formula  (\ref{btp}) holds at least  for every regular point $\xi$  of $\g^*$. Since   the set of regular elements of $\g^*$ is
 dense in  $\g^*$ (see  \cite[Proposition  19.7.5]{Yu}), the continuous map  $\xi\to[\diff_{\xi}F,\diff_{\xi}G]$ is zero at all points.
 This  establishes (\ref{btp}).
 The second point of the lemma follows directly from (\ref{fg09}).
\end{proof}
Let us now prove  Theorem \ref{tr1}.
\\

\begin{proof}\label{principal-proof}
%
{\em(1)}  Let  $K\in\FF(\g^*\times\g^*)$  and
 let   $(\xi,\eta)\in{\g^*\times\g^*}$. The  Poisson
$\RR$-bracket
  between  $\psi_1^* F=F\circ\psi _1$ and  $K$ at  $(\xi,\eta)$ is given by
\begin{eqnarray*}
{\pb{\psi ^*_1F ,K}}_{\RR}(\xi,\eta)
                       &=&\frac{1}{2}\can{(\xi,\eta),{[\RR\,\diff_{(\xi,\eta) }(F\circ\psi _1
                             ),\diff_{(\xi,\eta) }{K}]}}\\
                       &&+\frac{1}{2}\can{(\xi,\eta),{[\diff_{(\xi,\eta)}(F\circ\psi _1),\RR\,\diff_{(\xi,\eta) }{K}]}}.
\end{eqnarray*}%
According to (\ref{T15}) and   (\ref{TR}), we obtain $\diff_{(\xi,\eta)} (F\circ
\psi_1)=(\diff_{\xi-\eta } F,\diff_{\xi-\eta} F)$  and
   $\RR\,\diff_{(\xi,\eta) }(F\circ\psi _1  ) =c(\diff_{\xi-\eta }
F,\diff_{\xi-\eta }
F)$, so 
\begin{equation}\label{Rel1}
{\pb{\psi ^*_1F ,K}}_{\RR}(\xi,\eta)
                             =\frac{1}{2}\can{(\xi,\eta),[(\diff_{\xi-\eta} F
                             ,\diff_{\xi-\eta } F),\RR\,
    \diff_{(\xi,\eta) }{K}+c\diff_{(\xi,\eta) }{K}]}.
\end{equation}%
Using  Formula (\ref{formula1}), we have that  $\RR\,\diff_{(\xi,\eta)
}{K}+c\diff_{(\xi,\eta) }{K}$ is of the form  $(x,x)$ for some $x\in\g$, the explicit expression of which  is not needed in this proof. Then
 (\ref {Rel1}) becomes
\begin{eqnarray*}
{\pb{\psi ^*_1F ,K}}_{\RR}(\xi,\eta)
&=&\frac 12\can{(\xi,\eta),[(\diff_{\xi-\eta } F,\diff_{\xi-\eta } F),(x,x)]}\\
                            &=&\frac 12 \can{\xi-\eta,[\diff_{\xi-\eta }{F},x]}\\
                            &=&0,
\end{eqnarray*}
 where we used   that $F$ is an
$\Ad^*$-invariant function  on  $\g^*$ in the last line. This shows that $\psi ^*_1F$ is a Casimir on $((\gg)^*,{\PB}_{\RR})$.

{\em(2)}    Let  us prove that
the Poisson  $\RR$-bracket between
   $F\circ\psi _{\lambda}$ and  $G\circ\psi _{\gamma}$ at an arbitrary point  $(\xi,\eta)\in\g^*\times\g^*$ is equal to
    zero. According to (\ref{rlpb}),
\begin{equation}\label{Rel6}
\renewcommand{\arraystretch}{1.7}
\begin{array}{rcl}
{\pb{\psi _{\lambda}^*F,\psi _{\gamma}^*G}}_{\RR}(\xi,\eta ) &=&\frac{1}{2}\can{(\xi,\eta),[\RR\,\diff_{(\xi,\eta) }(F\circ\psi
    _{\lambda}),\diff_{(\xi,\eta) }(G\circ\psi _{\gamma})]}\\
                                                   &&+\frac{1}{2}\can{(\xi,\eta ),[\diff_{
  (\xi,\eta)  }(F\circ\psi
    _{\lambda}),\RR\,\diff_{(\xi,\eta) }(G\circ\psi _{\gamma})]}.
\end{array}
\end{equation}%
Using  Formulae (\ref{T15}) and  (\ref{TR}), the
Poisson $\RR$-bracket  ${\pb{\psi _{\lambda}^*F,\psi _{\gamma}^*G}}_{\RR
    }(\xi,\eta)$ of   (\ref{Rel6}) becomes
   \begin{eqnarray*}
   {\pb{\psi _{\lambda}^*F,\psi _{\gamma}^*G}}_{\RR}(\xi,\eta)&=&
\frac{1}{2}\can{(\xi,\eta),[(\lambda -1)(R\,\diff_{\lambda\xi-\eta } F
    ,R\,\diff_{\lambda\xi-\eta } F ),(\gamma \diff_{\gamma \xi-\eta }
    G,\diff_{\gamma \xi-\eta } G)]}\\
&&+\frac{1}{2}\can{(\xi,\eta),[(c\diff_{\lambda \xi-\eta }
  F,c\lambda\diff_{\lambda \xi-\eta } F),(\gamma \diff_{\gamma
    \xi-\eta } G,\diff_{\gamma \xi-\eta } G)]}\\
     &&-((F,\l)\leftrightarrow (G,\gamma))\\
&=&
\frac{(\lambda -1)}{2}\can{\gamma\xi-\eta,[R\,\diff_{\lambda \xi -\eta }
    F,\diff_{\gamma \xi+\eta } G]}\\
&&+\frac{c\gamma}{2} \can{\xi,[\diff_{\lambda \xi-\eta } F,\diff_{\gamma \xi -\eta }
      G]}- \frac{c\lambda}{2}\can{\eta, [\diff_{\lambda \xi-\eta }
    F,\diff_{\gamma \xi-\eta } G ]}\\
 &&-((F,\l)\leftrightarrow (G,\gamma))\\
&=&\frac{c\gamma}{2} \can{\xi,[\diff_{\lambda \xi-\eta } F,\diff_{\gamma \xi -\eta }
      G]}- \frac{c\lambda}{2}\can{\eta, [\diff_{\lambda \xi-\eta }
    F,\diff_{\gamma \xi-\eta } G ]}\\
 &&-((F,\l)\leftrightarrow (G,\gamma)),\\
\end{eqnarray*}%
 where we have used the
$\Ad^*$-invariance of $F$ and $G$  on    $\g^*$ to simplify the expression. According to  Formula (\ref{T18}), all four terms of the previous expression vanish if $\l\neq\gamma$. If  $\l=\gamma$, the bracket  ${\pb{\psi _{\lambda}^*F,\psi
    _{\gamma}^*G}}_{\RR}$  becomes
\begin{eqn*}
 {\pb{\psi _{\l}^*F,\psi _{\l}^*G}}_{\RR}(\xi,\eta)
={c\l} \can{\xi-\eta,[\diff_{\lambda \xi-\eta } F,\diff_{\l \xi -\eta }
      G]},
\end{eqn*}%
which is zero in view of Formula  (\ref{btp}).

We now suppose  that $F$ and $G$ are  polynomials of degrees $l$ and $k$, so for all
 $(\xi,\eta)\in\g\*\times\g^*$, $F\circ \psi_{\l}(\xi,\eta)=\sum_{i=0}^{l}\l^i F_i(\xi,\eta)$
 and $G\circ \psi_{\gamma}(\xi,\eta)=\sum_{j=0}^{k}\gamma^jG_j(\xi,\eta)$; hence
\begin{equation}\label{eqsumFG}
{\pb{F\circ\psi_\l,G\circ\psi_{\gamma}}}_{\RR}=\sum_{i=0}^{l}\sum_{j=0}^{k}\l^i\gamma^j{\pb{F_i,G_j}}_{\RR}.
\end{equation}%
As  we  just showed, the term on   the left hand side  of  Equation
(\ref{eqsumFG}) is  zero, for all $\l,\gamma\in\F$. Therefore, all the  coefficients
  of
 $\sum_{i=0}^{l}\sum_{j=0}^{k}\l^i\gamma^j{\pb{F_i,G_j}}_{\RR} $ are  zero.

{\em(3)} We  assume that   $c=1$. Let    $F,G\in\FF(\g^*)$ and let
  $(\xi,\eta)$ be a  point in   $\g^*\times\g^*$. From
 (\ref{Rel6}), (\ref{TR}) and  (\ref{T15}), it follows that the Poisson
$\RR$-bracket between  $F\circ\psi_{1}$ and  $G\circ\psi_{1}$ at
$(\xi,\eta)$ is given by
\begin{eqn*}
{\pb{F\circ\psi_1,G\circ\psi_1}}_{\RR}(\xi,\eta )=\can{\xi-\eta,[\diff_{\xi-\eta } F,\diff_{\xi-\eta } G]}=\pb{F,G}(\psi _1(\xi,\eta )).
\end{eqn*}%
Hence $\psi_1^*$ is a Poisson map.

{\em(4)} Let  $K$ be  a  function on   ${\g}^*\times\g^*$  and
$(\xi,\eta)$  be an element of
  $\g^*\times\g^*$. A direct computation gives
\begin{eqnarray}\label{Rel12}
\X_{H\circ\psi _{\lambda}}(\xi,\eta)[K]&=&{\pb{K,H\circ\psi _{\lambda}}}_{\RR}(\xi,\eta )\nonumber\\
                                &=&\frac{1}{2}\can{(\xi,\eta
  ),[\RR\,\diff_{(\xi,\eta) } K,\diff_{(\xi,\eta) } (H\circ \psi _{\lambda})]+[\diff_{(\xi,\eta)
    } K,\RR\,\diff_{(\xi,\eta) } (H\circ \psi _{\lambda})]}.\nonumber\\
&&\label{Rel12}
\end{eqnarray}%
 In order to rewrite this formula, we  temporarily use  $(x,y):=\diff_{(\xi,\eta)}K$. According to  (\ref{formula1}),  $\RR\,\diff_{(\xi,\eta)}K=\RR(x,y)$ is given by
\begin{equation}\label{g1}
   \RR\,\diff_{(\xi,\eta)}K=\left(R(x-y),R(x-y)\right)+c(y,x).
\end{equation}%
 Formulae  (\ref{g1}) and  (\ref{T15}) allow one  to rewrite the first term of (\ref {Rel12}) in the  following manner:
\begin{eqnarray}\label{Rel13}
\lefteqn{\can{(\xi,\eta),[\RR \diff_{(\xi,\eta) } K,\diff_{(\xi,\eta)
      } (H\circ \psi _{\lambda})]}}\nonumber\\
&&=
\can{(\xi,\eta ),[(R(x-y),R(x-y)),(\lambda \diff_{\lambda \xi-\eta }
    H,\diff_{\lambda \xi-\eta } H)]}\nonumber\\
&&~~+c\can{(\xi,\eta ),[(y,x),(\lambda \diff_{\lambda \xi-\eta }
    H,\diff_{\lambda \xi-\eta } H)]}\\
&&=
\can{\l\xi-\eta,[R(x-y),\diff_{\lambda\xi-\eta
} H]}+c\can{\l\xi,[y,\diff_{\lambda\xi-\eta } H]}-
c\can{\eta,[x,\diff_{\lambda\xi-\eta } H]}\nonumber\\&&=
c\can{\eta,[y,\diff_{\lambda\xi-\eta }
    H]}-c\can{\l\xi,[x,\diff_{\lambda\xi-\eta } H]},\nonumber
\end{eqnarray}%
where, in the second line, we used,  three times, the $\Ad^*$-invariance of $H$. Thus
\begin{equation}\label{Rel16}
\can{(\xi,\eta),[\RR\, \diff_{(\xi,\eta) } K,\diff_{(\xi,\eta) } (H\circ
    \psi _{\lambda})]}
=\can{(\xi,\eta),[(c\l\diff_{\lambda\xi-\eta }
    H,c\diff_{\lambda\xi-\eta }H),(x,y)]}.
\end{equation}
By replacing the first term of (\ref {Rel12}) by its  expression given in  (\ref
{Rel16}) and by using  Formula (\ref {TR}) to express the
 second term of (\ref {Rel12}),   we rewrite  $\X
 _ {H\circ\psi _ {\lambda}} (\xi, \eta) [F]$
as follows:
\begin{eqnarray*}
\X_{H\circ\psi
  _{\lambda}}(\xi,\eta)[K]&=&\frac{1}{2}\can{(\xi,\eta),[(c\lambda
    \diff_{\lambda \xi-\eta } H,c\diff_{\lambda \xi-\eta } H),(x,y)]}\\
&&-\frac{(\lambda -1)}{2}\can{(\xi,\eta),[(R\,\diff_{\lambda \xi-\eta }
    H,R\,\diff_{\lambda \xi-\eta } H),(x,y)]}\\
&&-\frac{1}{2}\can{(\xi,\eta),[(c\diff_{\lambda \xi-\eta } H,c\lambda
    \diff_{\lambda \xi-\eta } H),(x,y)]}\\
&=&\frac{(\lambda -1 )}{2}\can{(\xi,\eta),[((cI-R)\diff_{\lambda
      \xi-\eta } H,-(R + c I)\diff_{\lambda \xi-\eta } H),\diff_{(\xi,\eta) } K]}.
\end{eqnarray*}
We then deduce Formula (\ref{champh}). When $\g$ is a Lie algebra splitting and $R=P_+-P_-$,  we have $c=1$, $R-cI=-2P_-$ and $R+cI=2P_+$. This implies that (\ref{champh}) gives
(\ref{relchampspli}).
\end{proof}
\section{The Liouville integrability of  the $2$-Toda lattice}\label{2toda}

In this section we define the $2$-Toda lattice for every complex simple Lie algebra and we prove   its Liouville  integrability.
  
 In order to  define  the $2$-Toda lattice for  every simple
Lie  algebra,  we need some notation. Let $\g$ be  a simple Lie algebra of
rank $\ell$, with Killing form  $\INN$. We choose  $\h$, a
Cartan subalgebra with roots system $\Phi$, and
$\Pi=(\a_1,\dots,\a_{\ell})$,  a system of simple roots with
respect to $\h$. For every $\a$ in $\Phi\backslash \{-\Pi,\Pi\}$, we denote by  $e_{\a}$  a
non-zero eigenvector associated with eigenvalue $\a$, and, for every
$1\leqslant i\leqslant \ell$,  we denote  by $e_i$ and $e_{-i}$  non-zero  eigenvectors
associated respectively  with  $\a_i$ and $-\a_i$.  The Lie algebra $\g=\sum_{k\in\Z}\g_k$  is  endowed with the natural
 grading (i.e., for every $k,l\in\Z$, $[\g_k,\g_l]\subset\g_{k+l}$)
  defined by $\g_0:=\h$ and,  for every  $k\in\Z$,  
$\g_k:=\langle e_{\a}\mid   \a\in\Phi,  |\a|=k\rangle $, for  $|\a|$  the  length of the root $\a$, i.e.,  $|\a|$ is $\sum_{i=1}^{\ell}a_i$ for   $\a=\sum_{i=1}^{\ell}a_i\a_i$.
In the sequel, we shall  use the following property: $\inn{\g_k}{\g_l}=0$ if $k+l\neq 0$.
 We introduce the following notation
\begin{eqn*}
\g_{<k}:=\sum_{i<k}^{} \g _i,&\qquad \g_{\leqslant k}:=\sum_{i\leqslant k}^{}\g _i,\\
\g_{>k}:=\sum_{i>k}^{} \g _i,&\qquad \g_{\geqslant k}:=\sum_{i\geqslant k}^{} \g_i.
\end{eqn*}%
Also, $\gp:=\g_{\geqslant 0}$ and  $\gm:=\g_{< 0}$.
\subsection{Definition of the $2$-Toda lattice}
 
The next definition gives again  the  definition given in (\ref{k1}) 
 when specialized to the case $\g=\Liesl{n}(\C)$, taking  for  $\h$ the  Lie  subalgebra of  diagonal matrices.
 \begin{definition}
 The \emph{$2$-Toda lattice associated with a simple Lie algebra  $\g$}
 is the  system of differential equations given by  the following Lax equations:
\begin{eqnarray}
\renewcommand{\arraystretch}{2.5}
\ds\pp{(L,M)}{t}&=&[(L_+,L_+),(L,M)],\label{LL19}\\
\ds\pp{(L,M)}{s}&=&[(M_-,M_-),(L,M)],\label{LL199}
\end{eqnarray}
where   $(L,M)$ is an element of the phase space of the $2$-Toda lattice  $\TP:=\g_{\leqslant 0}\times \g_{\geqslant -1}+(\sum_{i=1}^{\ell}e_i,0)$, where  $L_+:=P_+(L)$,  $M_-:=P_-(M)$, and  where $P_{\pm}$ is   the projection of $\g$ on $\g_{\pm}$.
\end{definition}
\subsection{The $2$-Toda lattice is a Hamiltonian system}\label{sect-hamiltonian}

We recall from Example \ref{splitting-gg} that when $\g=\gp\oplus\gm$ is
a  Lie algebra splitting  then  $\gg=\ggp\oplus\ggm$ is also a Lie algebra splitting, where
\begin{equation}
\ggp:=\{(x,x) \mid x\in\g \} ~~\textrm { et } ~~ \ggm:=\{(x,y) \mid
x\in\gm\textrm { et } y\in\gp\}.
\end{equation}%
Also,  for  every $(x, y)\in\gg$ we have $(x,y)=(x,y)_++(x,y)_-$, where
\begin{equation}\label{equ-decom}
\left\{
 \begin{array}{ccc}
(x,y)_+&=&(x_++y_-,x_++y_-)\in\ggp,\\
(x,y)_-&=&(x_--y_-,y_+-x_+)\in\ggm.
\end{array}
\right.
\end{equation}%
 Let $\RR$ be  the difference of the projections of $\gg$ on the Lie
 subalgebras $\ggp$ and $\ggm$. According  to  Example \ref{splitting-gg}, we have
 \begin{eqnarray}\label{rmdt2}
\RR(x,y)&=&(x_+-x_-+2y_-,y_--y_++2x_+)\label{rmdt2}\\
        &=&(R(x-y)+y,R(x-y)+x)\label{R-theorem},
\end{eqnarray}%
where $ R$ is the difference of the  projections  of $\g$ on $\gp$ and on  $\gm$.
We provide $\gg$ with the  following $\Ad$-invariant, non-degenerate symmetric
bilinear  form:
\begin{equation}\label{blform}
\begin{array}{ccccc}
{\INN}_2&:&\gg\times\gg         &\to    &\C\\
    & &((x_1,y_1),(x_2,y_2))&\mapsto&\inn{x_1}{x_2}-\inn{y_1}{y_2}.
\end{array}
\end{equation}%
We use it to  identify $\gg$ with its dual and we obtain according to (\ref{pdg2m}) a linear
Poisson structure on $\gg$, defined for every $F,G\in\FF(\gg)$ at  $(x,y)\in\gg$ by
\begin{equation}\label{mlpfd21}
{\pb{F,G}}_{\RR}(x,y)=\frac{1}{2}{\inn{(x,y)}{{[\RR\nabla_{(x,y)}{F},\nabla_{(x,y)}{G}]}+{[\nabla_{(x,y)}{F},\RR\nabla_{(x,y)}{G}]}}}_2,
\end{equation}%
 where  $\nabla_{(x,y)}{F}$ is  the  gradient of  $F$ at  $(x,y)\in\gg$
 (with respect to $\INN_2$), i.e., 
 \begin{eqn*}
{\inn{\nabla _{(x,y)}{F}}{(z,s)}}_2=\can{\diff _{(x,y)}{F},(z,s)}, \qquad \forall (z,s)\in\gg.
\end{eqn*}%

We show that the phase space $\TP$  is equipped with a Poisson structure
 and the equations of motion of the $2$-Toda lattice  are Hamiltonian.
\begin{proposition}\label{K6}
$\TP$ is a   Poisson submanifold of
$(\gg,{\PB}_{\RR})$.
\end{proposition}
We use the following lemma to show the above proposition.
\begin{lemma}\label{affine-submanifold}
Let $\g$ be  a Lie algebra equipped with a non-degenerate, symmetric,
 bilinear form $\INN$. Let $a\in\g$ and let $E$ be  a subspace of $\g$. We suppose that:\\
\emph{(1)} The orthogonal\footnote{Here the
orthogonality is with respect to  the form $\INN$.} $E^{\perp}$ of $E $ is a  Lie ideal  of $(\g,\LB)$.\\
\emph{(2)} For every $x,y\in\g$, we have $\inn{a}{[x,y]}=0.$\\
Then $a+E$ and $E$ are  Poisson submanifolds of $(\g,{\PB})$, equipped with its linear Poisson structure $\PB$.
 The map $x\to a+x$ is a Poisson isomorphism
 between them. Moreover   $E$ is Poisson isomorphic to   the Lie-Poisson manifold  $(\g/E^{\perp})^*$, 
also equipped with its linear Poisson structure.
\end{lemma}

\begin{proof}
We start to show that $E$ is a Poisson submanifold of $(\g,\PB)$.
 Let $\J:=\langle F\in\FF(\g) \mid F\equiv 0 \textrm{ on } E \rangle$.  To show that  $E$ is a Poisson submanifold 
 is equivalent to show that   $\J$  is Poisson  ideal. Let $F\in\J$ and let $x\in E$, we notice that 
$\nabla_x F \in E^{\perp}$. 
According to the first condition of the proposition, for every $G\in\FF(\g)$, we have 
 $\lb{\nabla_xf,\nabla_xG}\in E^{\perp}$. This implies that  $\pb{F,G}(x)=\inn{x}{ \lb{\nabla_xf,\nabla_xG}}=0$.
 Then $\J$ is a Poisson ideal.
   
We now show that $E+a$ is a Poisson submanifold. Let   $\I:=\langle F\in\FF(\gg) \mid F\equiv 0 \textrm{ on } E+a\rangle $, 
  $F\in\I$ and 
  let $x+a\in E+a$. Notice that  $\nabla_{x+a}F$ is an element of  $E^{\perp}$, which is according to first  condition of 
these   
proposition is a Lie ideal of $(\g,\LB)$. Then, for every  $G\in\FF(\g)$,  
$[\nabla_{x+a}F,\nabla_{x+a}G]\in E^{\perp}$. Let now  compute the Poisson bracket between  $F$ 
 and  $G\in\FF(\g^)$ at  a point  $x+a$.
\begin{eqnarray*}
\pb{F,G}(x+a)&=&\inn{x+a}{[\nabla_{x+a}F,\nabla_{x+a}G]}\\
                &=&\inn{x}{[\nabla_{x+a}F,\nabla_{x+a}G]}\\
                &=&0,
\end{eqnarray*}%
where we have used the second condition  of  proposition to justify the transition  from first to second  line
 and we have zero in the last line  because  $x\in E$  and  $[\nabla_{x+a}F,\nabla_{x+a}G]\in E^{\perp}$.

Let now show   the  translation  by  $a$ is a Poisson isomorphism between  $\g$  and  $\g$. We denote by  
  $T_{a}$ the  translation  by  $a$, defined for every   $x\in\g$, by 
 $T_{a}(x)=x+a$.  For every  $F,G \in \FF(\g)$  and every   $x \in \g$,  we have
 \begin{eqnarray*}  
\pb{F\circ T_{a},G \circ T_{a}}(x)&=&\inn{x}{[\nabla _{x} (F \circ
    T_{a}) , \nabla  _{x}(G \circ T_{a})]}\\
 &=& \inn{x}{[\nabla_{T_{a} (x)}F,\nabla_{T_{a} (x)}G]} \\
&=& \inn{ x +a}{ [\nabla_{T_{a} (x)}F,\nabla_{T_{a} (x)}G] }    \\
&=&   \pb{F,G}(T_{a} (x)), 
 \end{eqnarray*}%
where we have used Condition  (2) to justify the transition from second to third  line. Since  $E$
 is a Poisson submanifold of 
 $(\g,\PB)$  the  restriction  of the Poisson isomorphism  $T_{a}$ to 
 $E$  is also a Poisson  isomorphism on  its image. 
\end{proof}%
We now prove  Proposition \ref{K6}.\\
\begin{proof}
 Its easy to verify  that $\TP$ admits the following description as an affine subspace of $\gg$:
\begin{equation}
\TP=(e,e)+\Delta(\g_{0}\oplus \g_{-1})\oplus  \gg_{-},
\end{equation}%
where  $\Delta( \g_{0}\oplus \g_{-1}):=\{(x,x) \mid  x\in  \g_{0}\oplus \g_{-1}\}$ and $e=\sum_{i=1}^{\ell}e_i$.
We check  that   the two   assumptions of  Lemma \ref{affine-submanifold} are satisfied, with  $a=(e,e)$,
 $E=\Delta(\g_{0}\oplus \g_{-1})\oplus  \gg_{-}$, $\g=\gg$ and $\LB={\LB}_{\RR}$.
\\
\emph{(1)}  It is clear that
$E=\g_{\leqslant 0}\times \g_{\geqslant -1}$. The orthogonal of $E$ is  $E^{\perp}=\g_{<0}\times\g_{>1}$,
  which  is a subspace of  $ \gg_-$ and  it is an  ideal of    $(\gg,\LB_\RR)$, because, for every
   $(x,y)\in\g_{<0}\times\g_{>1}$
and  $(z,s)\in\gg$,
\begin{eqn*}
  {[(x,y),(z,s)]}_{\RR}\in[\g_{<0}\times\g_{>1},\g_{<0}\times \g_{\geqslant 0}]\subset\g_{< -1}\times\g_{>1}\subset \g_{< 0}\times\g_{>1}.
\end{eqn*}%
\emph{(2)} For every   $(x,y),(x',y')\in\gg$, we have
\begin{equation}\label{chanbase}
\inn{(e,e)}{{[(x,y),(x',y')]}_{\RR}}={\inn{(e,e)}{{[{(x,y)}_+,{(x',y')}_+]}-{[{(x,y)}_-,{(x',y')}_-]}}}_2.
\end{equation}
Since $ [{(x,y)}_+,{(x',y')}_+]=(u,u)$ for some $u\in\g$ it follows that the first term on the right side of Equation
(\ref{chanbase}) is zero. Also    ${[{(x,y)}_-,{(x',y')}_-]}\in\g_{<-1}\times\g_{>0}$, 
which is orthogonal to $(e,e)\in\g_1\times\g_1$; hence the second term on the right hand  side of Equation (\ref{chanbase}) is zero.

\end{proof}

\begin{proposition}\label{hamil}

Let $H, {\tilde H}\in\FF(\gg)$  be defined  at  every point    $(x,y)$ of      $\gg$ by
\begin{equation}\label{hjk45}
H(x,y):=\frac{1}{2}\inn{x}{x}~~ \textrm{ and }~~ {\tilde{ H}}(x,y):=\frac 12\inn{y}{y}.
\end{equation}%
 The Hamiltonian  vector  field
  $\X_{H}:={\Pb{H}}_{\RR}$  (resp. $\X_{{\tilde H}}:={\Pb{\tilde H}}_{\RR}$) is tangent to $\TP$ and
describes on $\TP$  the equation of  motion (\ref{LL19})
(resp. (\ref{LL199})) for  the  $2$-Toda lattice.
\end{proposition}
\begin{proof}
Let $F\in\FF(\gg)$ and let $(L,M)\in\TP$. We have
\begin{eqnarray}
\X_H[F](L,M)&=&{\pb{F,H}}_{\RR}(L,M)\nonumber\\
           &=&\frac 12{\inn{(L,M)}{\lb{\RR\nabla_{(L,M)}F,\nabla_{(L,M)}H}}}_2\nonumber\\
           & &+\frac 12{\inn{(L,M)}{\lb{\nabla_{(L,M)}F,\RR\nabla_{(L,M)}H}}}_2\nonumber\\
          &=&-\frac 12{\inn{\lb{(L,M),\RR\nabla_{(L,M)}H }}{\nabla_{(L,M)}F}}_2,\label{equ-vector}
\end{eqnarray}%
where we have used the $\Ad$-invariance on   $\gg$ of $H$ to justify the transition from second to third equality.
We deduce from Equation (\ref{equ-vector}) that $\X_H(L,M)=\frac 12 \lb{\RR\nabla_{(L,M)}H ,(L,M)}$. Hence
according to  Formula   (\ref{equ-decom}), we have
\begin{eqn*}
\X_{H}(L,M)=[(\nabla_{(L,M) }{H})_+,(L,M)]=[(L,0)_+,(L,M)]=[(L_+,L_+),(L,M)].
\end{eqn*}%
To show that $\X_{\tilde H}$ describes  Equation (\ref{LL199})  it suffices to repeat the same reasoning of  $\X_H$.

Since $\TP$ is a Poisson submanifold of $(\gg,\PB)$, 
the Hamiltonian vector fields $\X_H$ and $\X_{\tilde H}$ are tangent to $\TP$.
%
\end{proof}

\subsection{The integrability of the $2$-Toda lattice}

According to \cite[Theorem 7.3.8]{Dixmier}, for  every  simple Lie algebra $\g$  of rank $\ell$,
there exist $\ell$ homogeneous, independent, $\Ad$-invariant polynomials $P_1,\dots,P_{\ell}$ which   generate 
the algebra of  
$\Ad$-invariant  polynomial functions on $\g$ and  which are    of degree, respectively,
$m_1+~1,\dots,m_{\ell}+1$, where $m_1,\dots,m_{\ell}$ are the exponents  of $\g$
 (we note that $m_1\leqslant\dots\leqslant m_{\ell}$).
\smallskip
Each  $P_i$ induces    $m_i+2$ functions $F_{j,i}\in\FF(\gg)$, as follows:
\begin{equation}\label{Cm8}
P _i(\lambda x-y)=\sum_{0\leqslant j\leqslant m_i+1}^{}(-1)^{m_i+1-j}\l^jF_{j,i}(x,y).
\end{equation}%
Every  function $F_{j,i}$ for  $1\leqslant i\leqslant \ell$ and  $0\leqslant j\leqslant
m_i+1$ is
 homogeneous  of  degree $j$ with respect  to its first variable and of degree $m_i+1-j$ with respect  to its second
 variable.
 \begin{notation}
We denote   by $\FF$ the family of functions   on  $\gg$  given by
\begin{equation}\label{ff}
 \FF:=(F_{j,i}, 1\leqslant i\leqslant \ell \textrm{ and } 0\leqslant j\leqslant m_i+1).
\end{equation}%
\end{notation}
\begin{remark}\label{remhamil}
\emph{(1)}  The functions $F_{0,1}$ and $F_{2,1}$ are the Hamiltonians of the $2$-Toda lattice introduced in Proposition \ref{hamil}.\\
\emph{(2)}  The functions $F_{0,i}$ and $F_{m_i+1,i}$ for $1\leqslant i\leqslant \ell$
 are $\Ad$-invariant functions on $\gg$. According  to the   Adler-Kostant-Symes theorem
\cite[Theorem 4.37]{Liv}  they are in
involution with respect to  the bracket  ${\PB}_{\RR}$. Also  they are
 independent  on  $\gg$.  They are therefore a good candidates for giving  the Liouville integrability of the $2$-Toda lattice.
 However   their
cardinal $2\ell$ is very small compared to $\dim\TP-\frac{1}{2}\Rk(\TP,{\PB}_{\RR})$.%
\end{remark}

 Since, as we  noticed
in Remark \ref{remhamil}, the Hamiltonians of the $2$-Toda lattice appear among  the functions composing $\FF$,
 the next theorem gives the Liouville integrability of the $2$-Toda lattice.
\begin{theorem}\label{principal-theorem}
The  triplet $(\TP,\FF_{\arrowvert\TP},{\PB}_{\RR})$ is  an integrable system.
\end{theorem}
\begin{proof}
 According  to the definition of  integrability in the sense of Liouville (see \cite[Definition 4.13]{Liv}), to prove Theorem
\ref{tr1} we must show that
\begin{enumerate}
\item[(1)] $\FF_{\arrowvert\TP} $ is  involutive for the Poisson $\RR$-bracket
  ${\PB}_{\RR}$.
\item[(2)] $\FF_{\arrowvert\TP}$ is independent.
\item[(3)] The cardinal of the restriction of $\FF$ to $\TP$ satisfies
\begin{equation}\label{vericard}
 \card\FF_{\arrowvert\TP}=\dim\TP-\frac{1}{2}\Rk(\TP,{\PB}_{\RR}).
 \end{equation}%
\end{enumerate}
The proofs of these  three points are given in respectively  Proposition \ref{kh6},   Proposition \ref{indtt} and Proposition \ref{K10},
which are given in   the next three subsections.
\end{proof}
\subsubsection{The restriction of $\FF$ to $\TP$  is  involutive for    ${\PB}_{\RR}$}

Unlike in the Toda lattice case, it is not
 the Adler-Kostant-Symes theorem that gives us the commutativity of 
the considered family of functions $\FF$, which is not formed only by the $\Ad$-invariant functions.
 This commutativity arises from the fact our phase space is a Poisson submanifold  of
 $(\gg,{\PB}_{\RR})$, and from the second item of Theorem \ref{tr1}. 
\begin{proposition}\label{kh6}
The family of functions  $\FF_{\arrowvert\TP}$ is involutive for the Poisson structure~${\PB}_{\RR}$.
\end{proposition}
\begin{proof}
Since the polynomials   $P_1,\dots,P_{\ell}$ are    $\Ad$-invariant
on  $\g$ (and hence $\Ad^*$-invariant upon  identifying  $\g$ with $\g^*$), according to the second  point of Theorem \ref{tr1}    the family  $\FF$ is  involutive on $(\gg,{\PB}_{\RR})$.
Furthermore,    $\TP$ being a Poisson submanifold  of $(\gg,{\PB}_{\RR})$ (see Proposition  \ref{K6}),  the restriction of $\FF$ to  $\TP$ is  involutive.
%
\end{proof}
\subsubsection{The restriction of $\FF$ to $\TP$  is  independent}

We use  an unpublished result  of Ra\"is \cite{Rais}, which  establishes  the independence
 of a large family of functions on $\gg$.  We  state  this result below and   the proof is  in \cite[Section 1]{ddg}.
\begin{theorem}\label{indvki}
Let  $P_1,\dots,P_{\ell}$   be a generating family of homogeneous polynomials of the algebra of  $\Ad$-invariant polynomial functions on   $\g$.
 Let  $e$ and  $h$ be  two elements of $\g$
 such that  $e$ is regular and  $[h,e]=2e$.

 For every  $F\in\FF(\g)$, and every  $y\in\g$, we denote  by  $\diff^k_yF$
 the differential of  order  $k$ of $F$  at $y$.
 Denote by  $V_{k,i}$, for every   $1\leqslant i\leqslant \ell$ and  $0\leqslant k\leqslant
 m_i$,  the element of  $\g$ defined by
\begin{equation}\label{dpvki}
\inn{V_{k,i}}{z}=\can { \diff^{k+1}_hP_i, (e^k,z)},\qquad \qquad \forall z\in\g,
\end{equation}%
where, for every $x\in\g$ and $k\in\N$, $x^k$ is  shorthand for $(x,\dots,x)$ ($k$ times).
\\
\emph{(1)} The family   $\FF_1:=(V_{k,i}, 1\leqslant i\leqslant \ell$ and $ 0\leqslant k\leqslant m_i)$
is  linearly   independent.\\
\emph{(2)}
The subspace  generated by    $\FF_1$ is the Lie  subalgebra formed by the sum of  the all  eigenspaces of $ \ad_h$
  associated with positive or zero eigenvalues.
\end{theorem}
 We first prove, using the  first point  of Theorem \ref{indvki}, the independence of
$\FF$ (which is a family $\FF(\gg)$)  at a well chosen  point $(e,h)\in\TP\cap\gg$. Afterwards,  using the second point of
   Theorem \ref{indvki}, we show that  the restriction
 of $\FF$ to  the phase space $\TP$ of the  $2$-Toda  lattice is  also an independent family of functions.
\begin{proposition}\label{rre}
Let    $h\in\h$ be  such that  $[h,e]=2e$.
\\
 \emph{(1)} The polynomial functions  $F_{0,1}(x,y),\dots,F_{0,\ell}(x,y)$ only  depend on the second
   variable~$y$. Their differentials  at the  point $(e,h)$ are  independent;\\
\emph{(2)} We denote by $\pp{}{x}$ the partial  differential with respect  to the variable $x$.
The   $\frac12(\dim\g+\ell)$ partial derivatives  $\pp{F_{j,i}}{x}$, $1\leqslant
   i\leqslant \ell $ and   $ 1\leqslant j\leqslant m_i+1$  are  independent at the point $(e,h)$;\\
 \emph{(3)}  The  family   $ (\diff_{(e,h)}{F_{j,i}}$, $1\leqslant i\leqslant \ell, 0\leqslant j\leqslant
   m_i+1)$ is    independent.
 \end{proposition}
\begin{proof}
 \emph{(1)}   For every  $1\leqslant i\leqslant \ell$ and every $(x,y)\in\gg$,
 $F_{0,i}(x,y)$ is the  term of  degree
 $0$ in  $\l$ of  $P_i(\l x-y)$,  that is   $P_i(y)$. Hence,  the function $F_{0,i}$,
for  $ 1\leqslant i\leqslant\ell$, is a homogeneous  polynomial of degree  $m_i+1$  which only depends on the second  variable
 $y$. Moreover, according   to  two theorems of  Kostant
 \cite[Theorem 9]{Bertram} and  \cite[Theorem  5.2]{kostant},
the differentials of  the polynomials  $P_1,\dots,P_{\ell}$
 are   independent  at every  regular   point of  $\g$. In particular, they are independent  at  the regular\footnote{The element  $h$ of  $\g$  is
   regular. Indeed, since it verifies   $[h,e]=2e$ with $e$ regular, it  belongs to a principal $\Liesl{}(2)$ triple (see
   \cite[Theorem 32.1.5]{Yu}).}  point  $h$. In conclusion  the  differentials of the functions   $F_{0,1},\dots,F_{0,\ell}$ are  independent  at  $(e,h)$.
   \\
\emph{(2)}    For $1\leqslant i\leqslant \ell$,
the  Taylor Formula  applied to the   polynomial $P_i$,  $1\leqslant i\leqslant\ell$,  at  $\lambda
x-y$, yields
\begin{equation}\label{diffp}
P_i(\lambda x-y)=\sum_{k=0}^{m_i+1}(-1)^{m_i+1-k}\frac{\lambda ^{k}}{k!}\can{\diff_y
^kP_i,x^k}.
\end{equation}%
By identifying   the coefficients of   Equations  (\ref{diffp}) and   (\ref{Cm8}),   we obtain, for every  $1\leqslant i\leqslant \ell$ and  $  1\leqslant k\leqslant  m_i+1$, the equality
\begin{equation}\label{relimp}
F_{k,i}(x,y)=\frac{1}{k!}\can{\diff_y ^kP_i,x^k}.
\end{equation}
By   differentiating $F_{k+1,i}$ with respect to the  variable $x$, we obtain
\begin{equation}\label{vkii}
\can{\pp{F_{k+1,i}}{x}(x,y),z}=\frac{1}{k!}\can{\diff_y ^{k+1}P_i,(x^k,z)},\qquad \forall
k=0,\dots, m_i,\textrm{ and }\forall z\in\g.
\end{equation}%
In particular, according to   Theorem  \ref{indvki}, when  $(x,y)=(e,h)$, Equation (\ref{vkii}) becomes
\begin{equation}\label{vki}
\can{\pp{F_{k+1,i}}{x}(e,h),z}=\frac{1}{k!}\can{\diff_h ^{k+1}P_i,(e^k,z)}=\frac{1}{k!}\inn{V_{k,i}}{z}, \qquad  \forall z\in\g,
\end{equation}%
for every   $1\leqslant i\leqslant \ell$ and  $  1\leqslant k\leqslant  m_i+1$.
According to Theorem \ref{indvki},
the  family $(V_{k,i}, 1\leqslant i\leqslant \ell \textrm{ and  } 0\leqslant k\leqslant\ell)$
is  independent. This  implies  the
independence of the  family
 $(\pp{F_{k,i}}{x}(e,h), 1\leqslant i\leqslant\ell\textrm{ and } 1\leqslant k\leqslant
m_i+1)$.
\\
\emph{(3)} Denote by   $M$ the  matrix
\begin{eqnarray*}
 M&=&(\diff _{(e,h)}F_{k,i}, ~~0\leqslant k\leqslant m_i+1 \textrm{ and  } 1\leqslant i\leqslant \ell )\\
  &=&\left(
     \begin{array}{cc}
     \pp{F_{0,1}}{x}(e,h)&\pp{F_{0,1}}{y}(e,h)\\
     \vdots              &\vdots               \\
      \pp{F_{0,\ell}}{x}(e,h)&\pp{F_{0,\ell}}{y}(e,h)\\
      \pp{F_{1,1}}{x}(e,h)&\pp{F_{1,1}}{y}(e,h)\\
     \vdots              &\vdots               \\
      \pp{F_{m_{\ell}+1,\ell}}{x}(e,h)&\pp{F_{m_{\ell}+1,\ell}}{y}(e,h)
      \end{array}
      \right).
\end{eqnarray*}
According  to item {\em(1)}, the form of the  matrix  $M$ is
$$
M=\left(
\begin{array}{cc}
0&A\\
B&*
\end{array}
\right),
$$
where    $A=(\pp{F_{0,i}}{y}(e,h), 1\leqslant i\leqslant\ell)$ and
  $B=(\pp{F_{k,i}}{x}(e,h), 1\leqslant k\leqslant m_i+1\textrm{ and } 1\leqslant i\leqslant \ell )$.
According   to {\em(1)} and  {\em(2)}, the  $\ell$  rows
of the    matrix $A$ and the  $\frac12(\dim\g+\ell)$ rows
of the   matrix   $B$ are   independent. The independence of  the  rows of   $M$ yields
 the independence of the  differentials  of the  family of  functions
$\FF$ at
$(e,h)\in\TP $.
\end{proof}
\begin{proposition}\label{indtt}
The  restriction of the  family of  functions  $\FF$ to the phase space
$\TP$ of the  $2$-Toda  lattice   is an independent family.
\end{proposition}
\begin{proof}
It suffices to show that the  restrictions on the tangent space $\T_{(e,h)} \TP $
of the differentials
$$(\diff_{(e,h)} F_{j,i},~~1\leqslant i\leqslant\ell\textrm{ et }0\leqslant j\leqslant m_i+1)$$
are   independent. Since  the  tangent space  of  $\TP$ is
 $\g_{\leqslant
  0}\times\g_{\geqslant -1}$ and since $(\diff_{(e,h)}F_{0,i},1\leqslant i\leqslant \ell)$ vanish on $\g_{\leqslant 0}\times\{0\}$, it suffices to  show that:
\\
\emph{(a)} The restriction to  the space  $\{0\} \times  \h$
of the linear  forms
$$(\diff_{(e,h)} F_{0,i},~~1\leqslant i\leqslant\ell~ )$$
 is  an independent family.\\
\emph{(b)}
The restriction  to the space   $  \g_{\leqslant 0} \times \{0\}$
of the linear independent forms
$$(\diff_{(e,h)} F_{j,i},~~1\leqslant i\leqslant\ell\textrm{ and }1\leqslant j\leqslant m_i+1)$$
is  an independent family.

According  to Definition (\ref{Cm8}), $F_{0,i}(x,y)=P_i(y)$ for every
$ 1\leqslant i\leqslant \ell$ and every $(x,y)\in\gg$. Then  $\pp{F_{0,i}}{x}(e,h)=0$ and
\begin{equation}
\can{\pp{F_{0,i}}{y}(e,h),z}=\inn{\nabla_hP_i}{z},\qquad\forall z\in\g.
\end{equation}%
Since   $\nabla_hP_i\in\h$  (because  $[h,\nabla_h P_i]=0$ and   $h\in\h$ is  regular),
    the  restriction  to   $\h$ of the   family
$(\pp{F_{0,1}}{y}(e,h),\dots, \pp{F_{0,\ell}}{y}(e,h))$ is  an  independent  family of linear  forms. This implies
 the independence of the   restriction to  $\{0\}\times \h$
of  $ \diff_{(e,h)} F_{0,1},\dots,\diff_{(e,h)} F_{0,\ell}$.

\smallskip

According to  equation (\ref{vki}),  for every  $1\leqslant i\leqslant
\ell$ and  $0\leqslant k\leqslant \ell$,
\begin{equation}\label{regr-forl}
\inn{V_{k,i}}{z}=k!\can{\pp{F_{k+1,i}}{x}(e,h),z}.
\end{equation}%
The second  point of Theorem \ref{indvki} shows that  the subspace generated by the family
  $(V_{k,i}$,
$1\leqslant i\leqslant \ell $ and  $0\leqslant k\leqslant m_i)$ is contained   in  the  Lie subalgebra  obtained by 
summing  the eigenspaces of  $\ad_h$, associated with the  nonnegative eigenvalues,
 which is  in our case   $\g_{\geqslant 0}$. This  proves, using Equation
(\ref{regr-forl}), that the  restriction  to  $\g_{\leqslant 0}$ of the family of linear forms
$(\pp{F_{k,i}}{x}(e,h)$,  $1\leqslant i\leqslant \ell$ and  $1\leqslant k\leqslant m_i+1)$
is an independent family. Therefore the  restriction to  $\g_{\leqslant 0}\times \{0\}$ of
$(\diff_{(e,h)}F_{k,i}, 1\leqslant i\leqslant \ell \textrm{ and  } 1\leqslant i\leqslant m_i+1)$ is an independent family.

\end{proof}
\subsubsection{The exact number of functions}

According to Equation (\ref{ff}),  the  cardinal of  $\FF$ is related   to the   exponents $m_i$,  $1\leqslant i\leqslant \ell$, as follows
 $$\card\FF=\sum_{i=1}^{\ell} (m_i+2)=\sum_{i=1}^{\ell}m_i+2\ell. $$
 Since $\sum_{i=1}^{\ell} m_i=\frac{1}{2}(\dim\g-\ell) $ (see \cite[Theorem 7.3.8]{Dixmier}) and  
$\FF_{\arrowvert\TP} $ is an independent family,   we have    $\card\FF_{\arrowvert\TP}= \frac{1}{2}
  (\dim\g+3\ell)$.
The dimension of $\TP$ is equal to  $\dim\g+2\ell$. In conclusion,  the  relation below is satisfied:
$$
\card\FF_{\arrowvert\TP}= \dim\TP-\frac{1}{2}\Rk(\TP, {\PB}_{\RR})
$$
if and only if  $\Rk(\TP,{\PB}_{\RR})=\dim\g+\ell$.  We need therefore to prove this last result, which will be done in Proposition  \ref{K10} below.
\paragraph{The rank of the restriction of ${\PB}_{\RR}$ to  $\TP$:}
The purpose of this part is to compute the rank of   the Poisson manifold  $(\TP,{\PB}_{\RR})$, i.e., the maximum of the rank at $x$ of the Poisson $\RR$-bracket for every $x\in\TP$. 
 We begin by establishing an isomorphism between the Poisson manifold $(\TP,{\PB}_{\RR})$  and  a product Poisson manifold.

Before doing this,  notice that it follows from Lemma
 \ref{affine-submanifold} applied  to $a=0$ and $E=\g_0\oplus\g_1$, $\g=\g$ and $\LB={\LB}_R$  that (the subspace)
 $\g_0\oplus\g_{1}$ is a Poisson submanifold of $(\g,\PB_R)$, since  the orthogonal of $\g_0\oplus\g_1$
 is  a Lie  ideal of $(\g,\LB_R)$.
\begin{proposition}\label{isomoTP}
 The Poisson manifold  $(\TP,{\PB}_{\RR})$ is  isomorphic to the product Poisson manifold
  $(\g^*,\PB )\times (\g_{0} \oplus \g_{1},{\PB}_R)$, where  $\PB$ is the Lie-Poisson bracket on  $\g^*$,
 $R$ is the  difference of projections on $\gp$ and $\gm$, and  ${\PB}_R$
 is the  Poisson $R$-bracket  on  $\g$ (restricted to  $\g_0\oplus\g_1$).
\end{proposition}
\begin{proof}
As we saw in the proof of  Proposition  \ref{K6},
 $a:=(e,e)$ and  $E:=\gg_-\oplus
\Delta(\g_0 \oplus \g_{-1})$ satisfy Conditions \emph{(1)} and \emph{(2)} of Lemma \ref{affine-submanifold};
hence $\TP=(e,e)+E$ and $E$ are isomorphic Poisson submanifolds  of $(\gg,{\PB}_{\RR})$. 
According to  Lemma \ref{affine-submanifold} again,
 $E$ is  Poisson isomorphic to
$(\gg/E^{\perp})^* $, where  $\gg/E^{\perp} $ is
 equipped with the quotient of the Lie bracket $\LB_{\RR}$ with respect to the Lie ideal
$E^{\perp}$. In conclusion, $(\TP,{\PB}_{\RR})$ is Poisson isomorphic to $((\gg/E^{\perp})^*,{\PB}') $, where ${\PB}'$ is the Lie Poisson 
bracket of $(\gg/E^{\perp})^* $.
Let us
 determine the    quotient Lie algebra  $\gg/E^{\perp} $.
The Lie algebra  $\gg$ (endowed  with  the bracket $\LB_\RR$) is  isomorphic to  the direct  sum
of the  Lie algebras   $\g $, $ \g_{\leqslant -1}$ and  $ \g_{\geqslant 0}$,
the isomorphism being  given, for every  $x  \in  \g, y_-  \in  \g_-, y_+ \in \g_+ $, by
  \begin{equation}\label{isomorphic}
  \begin{array}{ccc}
\g\times\gm\times\gp& \simeq &\gg\\
 (x,y_-,y_+) &\to& (x,x)+(y_-,y_+).
\end{array}
 \end{equation}%
 The Lie algebra  isomorphism (\ref{isomorphic})
 identifies the subspace   $E^{\perp}=\gm \times\g_{>1 }\subset  \gg $
 with  $(0,\gm,\g_{\geqslant 2})$. Hence  the  quotient $\gg/E^{\perp} $
is the  direct sum  of the  Lie algebra $\g$   (endowed with the usual bracket)
and  the quotient Lie  algebra  $\g_+/\g_{\geqslant 2} $.

We conclude that  $\TP $ is Poisson  isomorphic to the  dual of the direct sum
Lie algebra  $\g \oplus (\g_+/\g_{\geqslant 2}) $,
endowed with the Lie-Poisson  structure, and then is isomorphic to the product of the Poisson manifold
 $\g^*$ (endowed with the  Lie-Poisson structure)
with   $(\g_+/\g_{\geqslant 2})^*$ (endowed with  the  Lie-Poisson structure).
To complete the proof, it suffices to recall  (see Lemma  \ref{affine-submanifold}) that   $\g_0 \oplus
\g_1$ is a Poisson submanifold of  $(\g, \PB_R)$,
isomorphic to the Lie-Poisson  structure on the dual of  the quotient Lie algebra  $\g/(
\g_0\oplus \g_1)^{\perp} = \g_+/\g_{\geqslant 2}$.
\end{proof}
\begin{proposition}\label{K10}
The  rank of  the restriction of  the Poisson  $\RR$-bracket  ${\PB}_{\RR}$ to the manifold  $\TP$  is  $\dim\g+\ell$.
As a consequence, the following relation is satisfied:
$$
 \card\FF_{\arrowvert\TP}= \dim\TP-\frac{1}{2}\Rk(\TP,{\PB}_{\RR}).
$$ 
\end{proposition}

\begin{proof}
 According to Proposition  \ref{isomoTP} the Poisson submanifold  $(\TP,{\PB}_{\RR})$ is
  isomorphic  to the product  manifold    $(\g^*,\PB )\times(\g_0\oplus\g_1,{\PB}_R)$.
 This result  proves that the restriction of the  rank of Poisson $\RR$-bracket  ${\PB}_{\RR}$  to
 the manifold  $\TP$ is the sum  of the  rank of the  Lie-Poisson  structure
 on  $\g^*$, which is
  $\dim\g-\ell$ (see \cite[Proposition 29.3.2]{Yu}), and the   rank of the  Poisson $R$-bracket on
 $\g_0\oplus\g_1$.

 We calculate the latter  rank and show that it is $2\ell$, which finishes the proof.
Let  $(z_1,\dots,z_{2\ell})$ be the  coordinate system on
$\g_0\oplus\g_1$, defined by
$$
\left\{
\begin{array}{rcl}
z_i(x)&=&\inn{h_i}{x},\\
z_{i+\ell}(x)&=&\inn{e_{-i}}{x},
\end{array}
\qquad\qquad \forall i\in\{1\dots,\ell\},
\right.
$$%
where as before, for  $1\leqslant i\leqslant \ell$, the element  $e_{-i}$ is a  non-zero  eigenvector
 associated with the root  $-\a_i$.
The Lie-Poisson brackets  between these coordinate   functions are given by the following formulae:
$$
\left\{
\begin{array}{rcl}
{\pb{z_i,z_j}}_R&=&{\pb{z_{i+\ell},z_{j+\ell}}}_R=0,\\
{\pb{z_{i+\ell},z_j}}_R&=&c_{ij}z_{i+\ell},\\
\end{array}
\qquad\qquad\qquad \forall i,j\in\{1,\dots,\ell\},
\right.
$$%
where   $C:=(c_{ij})_{1\leqslant i,j\leqslant\ell}$ is the  Cartan   matrix of   $\g$.
Then the Poisson matrix
$M=({\pb{z_i,z_j}}_R)_{1\leqslant i,j\leqslant 2\ell}$ is equal to
$$
M=\left(
\begin{array}{cc}
0&-A^T\\
A&0
\end{array}
\right),
$$%
where  $A:=vC$ and  $v:=diag(z_{1+\ell},\dots,z_{2\ell})$.
Since  $C$ is  invertible,
the   rank of   $A$ at any  point for which  $(z_{1+\ell},\dots, z_{2\ell})=(1,\dots,1)$  is  $\ell$; therefore the
 rank of   $M$ is  $2\ell$, which implies that the  rank of the  Poisson $R$-structure
on
   $\g_0\oplus\g_1$ is  $2\ell$.
\end{proof}
\subsection{The integrability   for the quadratic Poisson $\RR$-bracket}
\label{sec:quadratic}

In this subsection we  study the integrability 
of the $2$-Toda lattice on  some Lie algebra $\g$  with respect to a
  quadratic Poisson $\RR$-bracket.  To construct a quadratic Poisson $\RR$-bracket on $\gg$ it is necessary that  
  $\g$ is an associative algebra of finite dimensional  and $\RR$ and its  antisymmetric part $\RR_-$ are solutions of (mCYBE).
 For this it suffices to choose $\g=\gl{n}(\C)$ and $\RR(x,y)=(R(x-y)+y,R(x-y)+x)$, where   where  $R=P_+-P_-$,  $P_+$
(resp. $P_-$) being the projection of $\gl{n}(\C)$ on the subalgebra of upper (resp. strictly lower)
triangular matrices 
(we will show again later that $\RR_-$ is a solution of (mCYBE) of $\gg$).  

 In this subsection we show the integrability
of a system of equations which is exactly the system (\ref{k1}),
up to the fact that we do no longer assume the matrices $L,M$ to be traceless. 
We denote by  ${\TP}'$ the phase space  of the $2$-Toda lattice, i.e.,  
 \begin{eqn*}
 \begin{array}{c}
{\TP}':= \{ (L,M)\in\gl{n}(\C)\times\gl{n}(\C)\mid \\(L,M)=
\left(
\begin{array}{c}
\left(
\begin{array}{cccc}
a_{11}&1     &         &0     \\
a_{21}&a_{22}&\ddots   &      \\
\vdots&      &\ddots   &1     \\
a_{n1}&\cdots&a_{n,n-1}&a_{nn}
\end{array}
\right)
,
\left(
\begin{array}{cccc}
b_{11}&\cdots&\cdots   &b_{1n}\\
b_{21}&\ddots&         &\vdots\\
      &\ddots&\ddots   &\vdots\\
0     &      &b_{n,n-1}&b_{nn}
\end{array}
\right)
\end{array}
\right) \}
\end{array}
\end{eqn*}%
and we name as the  \emph{$2$-Toda lattice on $\gl{n}(\C)$} the system of differential equations (\ref{k1})
 with the constraint $(L,M)\in{\TP}'$.

 Before studying the integrability of the $2$-Toda lattice with respect to the  quadratic Poisson 
$\RR$-bracket we will study the integrability of the latter system with respect to  the  linear 
Poisson $\RR$-bracket.
 \subsubsection{The Liouville integrability of the $2$-Toda lattice on $\gl{n}(\C)$ for the linear Poisson $\RR$-bracket}

We equip $\gl{n}(\C)\times\gl{n}(\C)$ with 
the $\Ad$-invariant, non-degenerate, symmetric, bilinear  form $\INN_2$, defined for every 
$(x,y), (x',y')\in\gl{n}(\C)\times\gl{n}(\C)$ by
\begin{equation}\label{form-square}
{\inn{(x,y)}{(x',y')}}_2:=\inn{x}{x'}-\inn{y}{y'}=\Trace(xx')-\Trace(yy').
\end{equation}%
As in Subsection \ref{sect-hamiltonian} we consider $\RR(x,y)=(R(x-y)+y,R(x-y)+x)$ 
and we consider the linear Poisson $\RR$-bracket\footnote{According to  Corollary \ref{pror}, 
the endomorphism $\RR$ is an $\RR$-matrix on $\gl{n}(\C)\times\gl{n}(\C)$.}, 
 defined for every $F,G\in\FF(\gl{n}(\C)\times\gl{n}(\C))$ 
at  $(x,y)\in\gl{n}(\C)\times\gl{n}(\C)$, by
$$
{\pb{F,G}}_{\RR}(x,y)=\frac 12({\inn{(x,y)}{[\RR \nabla_{(x,y)}F,\nabla_{(x,y)}G]}}_2
+{\inn{(x,y)}{[ \nabla_{(x,y)}F,\RR\nabla_{(x,y)}G]}}_2).
$$  
By using  the proofs of Propositions \ref{K6} and \ref{hamil}  we show that the phase space 
${\TP}'$ is a Poisson submanifold of $(\gl{n}(\C)\times\gl{n}(\C),\PB_{\RR})$ and $\X_H:={\Pb{H}}_{\RR}$ and 
$\X_{H'}:={\Pb{H'}}_{\RR}$ where 
$H(x,y)=\frac 12 \Trace x^2$ and $H'(x,y)=\frac 12\Trace y^2$ describes on ${\TP}'$ the equations of motion of
the $2$-Toda lattice.

Let us now study the integrability of  the $2$-Toda lattice on $\gl{n}(\C)$. Let $P_i$, 
for every $i\in\N$  be the $\Ad$-invariant 
function of $\gl{n}(\C)$ defined for all $x\in\gl{n}(\C)$
 by $P_i(x)=\frac{1}{i+1}\Trace(x^{i+1})$. We define $F_{j,i}\in\FF(\gl{n}(\C)\times \gl{n}(\C))$ by  
 $P_i(\l x-y)=\sum_{j=0}^{i+1}(-1)^{m_i+1-j}\l^jF_{j,i}(x,y)$ and we define
 $$
{\FF}'=(F_{j,i}, 0\leqslant i\leqslant n-1\textrm{ and } 0 \leqslant j \leqslant  i+1).
 $$ 
Note that the  functions that  make up  ${\FF}'$ are
the  functions that make up    $\FF$ together with 
the   functions $ F_{1,0}(x,y) := \Trace(x)$ and $F_{0,0}(x,y) := \Trace (y)$. We have the following proposition.
\begin{proposition}\label{linear-proof}
The triplet $({\TP}',\FF'_{|_{{\TP}'}}, \PB_\RR)$ is an integrable system.
\end{proposition} 
\begin{proof}
According to the item \emph{(2)} of   Theorem \ref{tr1} the family $\FF'$ is involutive for the Poisson bracket $\PB_{\RR}$. 
Since   ${\FF}'$ is an  independent family  on ${\TP}'$ 
(this follows from the independence of the differentials $\diff F_{0,0}$ and $\diff F_{1,0}$  at all points 
of $\gl{n}(\C)\times\gl{n}(\C)$ and the independence 
 of $\FF$ on  the submanifold
$\TP$, which is the submanifold of ${\TP}'$ defined by  $ F_{1,0}=F_{0,0}=0$), and 
since the cardinal of ${\FF}'$ is $\frac{n(n+3)}{2}$
we have, according to \cite[Proposition  4.12]{Liv}, the inequality 
\begin{equation}\label{eq:rqnksense2}
\Rk ({\TP}',{\PB}_{\RR})\leqslant 2({\rm dim}({\TP}')- {\rm Card}({\FF}'))= n^2+n-2.
\end{equation}%
Furthermore  $\Rk({\TP}',{\PB}_{\RR})\geqslant \Rk({\TP},{\PB}_{\RR})=n^2+n-2$. 
This implies that 
 the rank of the restriction
 to ${\TP}'$ of $\PB_\RR^Q $ is exactly $n^2+n-2$. Then 
 $ \card{{\FF}'}_{\arrowvert{\TP}'} =\dim{\TP}'-\frac{1}{2}\Rk({\TP}',{\PB^Q}_{\RR})$. This completes the proof.
 \end{proof}
\subsubsection{The integrability of the $2$-Toda lattice on $\gl{n}(\C)$ for the quadratic Poisson $\RR$-bracket}

Let the quadratic $\RR$-bracket  be defined for  every $F,G\in\FF(\gl{n}(\C)\times\gl{n}(\C))$ at 
 $(x,y)\in\gl{n}(\C)\times\gl{n}(\C)$, by
\begin{eqnarray}\label{quadrPB}
{\pb{F,G}}^{Q}_{\RR}(x,y)&:=&\frac{1}{2}{\inn{[(x,y),\nabla_{(x,y)}F]}{\RR((x,y)\nabla_{(x,y)}G+\nabla_{(x,y)}G(x,y))}}_2\nonumber\\
                          &&-(F\leftrightarrow G).
\end{eqnarray}%
Since $\RR_-(x,y)=\frac 12(\RR-\RR^*)(x,y)=(R_-(x-y)+y,R_-(x-y)+x)$ and 
 since both the endomorphism $R=P_+-P_-$ and its  antisymmetric part $R_-=\frac{1}{2}(R-R^*)$
 are solutions of (mCYBE) of $\gl{n}(\C)$ of constant $c=1$,
it follows from Corollary \ref{pror} that $\RR$ and $\RR_-$ are solutions of (mCYBE) of $\gl{n}(\C)\times\gl{n}(\C)$ 
of constant $c=1$. According to \cite[Section 4]{Parmentier}
the $\RR$-bracket  (\ref{quadrPB}) is indeed  a Poisson bracket, that we call \emph{the  quadratic Poisson 
$\RR$-bracket (on $\gl{n}(\C)\times\gl{n}(\C)$}).

By a direct computation in coordinates we can  prove the following proposition.
\begin{proposition}\label{prop:quadra1}
The phase space ${\TP}'$ of the $2$-Toda lattice on $\gl{n}(\C)$  is a Poisson submanifold of 
$(\gl{n}(\C)\times\gl{n}(\C),{\PB}_{\RR}^Q)$.
\end{proposition}
We show the following proposition.
\begin{proposition}\label{prop:quadra2}
Let $P_i$, for every $i\in\N$  be the $\Ad$-invariant functions of $\gl{n}(\C)$ defined for all $x\in\gl{n}(\C)$
by $P_i(x)=\frac{1}{i+1}\Trace(x^{i+1})$ and let $\phi_{\l}:\gl{n}(\C)\times\gl{n}(\C)\to\gl{n}(\C)$, for every $\l\in\F$, be 
   defined
by $\phi_{\l}:(x,y)\to\l x-y$. \\
\emph{(1)} For every $i,j\in\N$ and every $\l, \gamma\in\C$, the functions $P_i\circ\phi_{\l}$ and $P_j\circ\phi_{\gamma}$
are in involution for $\PB_\RR^Q$.\\
\emph{(2)} The Hamiltonian vector field $\X_{P_i\circ\phi_{\l}}^Q:={\Pb{P_i\circ\phi_{\l}}}_{\RR}^Q$ is given by
 \begin{equation}\label{vectorquad}
\X_{P_i\circ\phi_{\l}}^Q=-[(x,y),((R-I)(\l x-y)^{i+1},(R+I)(\l x-y)^{i+1})].
\end{equation}%
\end{proposition}
\begin{proof}
 \emph{(1)} According to (\ref{quadrPB}) and (\ref{T15}), for every  $x,y\in\g$,
\begin{eqnarray*}
\lefteqn{  {\pb{P_i\circ\phi_{\lambda},P_j\circ\phi_{\mu}}}^Q_{\RR}(x,y)}\\
&=&\frac 12 {\inn{[(x,y),(\lambda(\lambda x-y)^i,(\lambda x-y)^i)]}{\RR(\mu x(\mu x-y)^j, y(\mu x-y)^j)}}_2\\
&&+\frac 12 {\inn{[(x,y),(\lambda(\lambda x-y)^i,(\lambda x-y)^i)]}{\RR(\mu(\mu x-y)^jx, (\mu x-y)^jy)}}_2\\
&&-(j,\mu)\longleftrightarrow(i,\lambda).
\end{eqnarray*}%
By replacing $\RR$ by its expression $\RR(x,y)=(R(x-y)+y,R(x-y)+x)$, we obtain
\begin{eqnarray*}
\lefteqn{ {\pb{P_i\circ\phi_{\lambda},P_j\circ\phi_{\mu}}}^Q_{\RR}(x,y)}\\
&=& {\inn{[(x,y),(\lambda(\lambda x-y)^i,(\lambda x-y)^i)]}{(R((\mu x-y)^{j+1}),R((\mu x-y)^{j+1}))}}_2\\
&&+ \frac{1}{2} {\inn{[(x,y),(\lambda(\lambda x-y)^i,(\lambda x-\mu)^i)]}{( y(\mu x-y)^j,\mu x(\mu x-y)^j)}}_2\\
&&+ \frac{1}{2} {\inn{[(x,y),(\lambda(\lambda x-y)^i,(\lambda x-\mu)^i)]}{( (\mu x-y)^jy,\mu(\mu x-y)^jx)}}_2\\
&&-(j,\mu)\longleftrightarrow(i,\lambda).
\end{eqnarray*}
By using  Formula (\ref{form-square}) we obtain 
\begin{eqnarray*}%
 {\pb{P_i\circ\phi_{\lambda},P_j\circ\phi_{\mu}}}^Q_{\RR}(x,y)
&=&\inn{[\l x-y,(\l x-y)^i]}{R((\mu x-y)^{j+1})}\\
&&+\frac {1}{2}  \inn{[\lambda x,(\lambda x-y)^i]}{ y(\mu x-y)^j+(\mu x-y)^jy}\\
&&-\frac{1}{2} \inn{[y, (\lambda x-y)^i]}{\mu x(\mu x-y)^j+\mu(\mu x-y)^jx)} \\
&& - (j,\mu)\longleftrightarrow(i,\lambda).
\end{eqnarray*}%
Since $[\l x,(\l x-y)^i]=[y,(\l x-y)^i]$,
\begin{eqnarray*}
 {\pb{P_i\circ\phi_{\lambda},P_j\circ\phi_{\mu}}}^Q_{\RR}(x,y)
&=&-\frac{1}{2} \inn{[y, (\lambda x-y)^i]}{2(\mu x-y)^{j+1}} - (j,\mu)\longleftrightarrow(i,\lambda)\\
&=&\inn{y}{[ (\lambda x-y)^i,(\mu x-y)^{j+1}]} -  (j,\mu)\longleftrightarrow(i,\lambda)\\
&=&0,
\end{eqnarray*}
where we used  Lemma \ref{commute1} to provide the last line.\\
\emph{(2)} Let $K$ be a function of $\gg$; we denote by $(a,b)=\nabla_{(x,y)}K$, according to  (\ref{quadrPB})
and~(\ref{T15}),
\begin{eqnarray}
\lefteqn{\X_{P_i\circ\phi_{\l}}^{Q}(x,y)[K]}\nonumber \\
&=&\frac 12\inn{[(x,y),(a,b)]}{\RR(\l x({\l x-y})^i+\l(\l x-y)^ix,y(\l x-y)^i+
(\l x-y)^i y)}\nonumber\\
                           &&-\frac 12\inn{[(x,y),(\l(\l x-y)^i,(\l x-y)^i)]}{\RR(xa+ax,yb+by)}\label{qvf}.
\end{eqnarray}
We introduce the shorthand 
$B=xa+ax-yb-by$. 
By using the expression for $\RR$  the Hamiltonian vector field (\ref{qvf})  becomes
\begin{eqnarray*}
 \X_{P_i\circ\phi_\l}^{Q}(x,y)[K]
&=& {\inn{[(x,y),(a,b)]}{(R((\l x-y)^{i+1}),R((\l x-y)^{i+1}))}}_2\\
 &&+\frac{1}{2}{\inn{[(x,y),(a,b)]}{(y(\l x-y)^i+(\l x-y)^iy,\l x(\l x-y)^i+\l(\l x-y)^ix)}}_2\\
 &&-\frac 12 {\inn{[(x,y),(\l (\l x-y)^i,(\l x-y)^i)]}{(R(B),R(B))}}_2\\
 &&-\frac{1}{2}{\inn{[(x,y),(\l(\ x-y)^i,(\l x-y)^i)]}{(yb+by,xa+ax)}}_2\\
 &=&- {\inn{[(x,y),(R((\l x-y)^{i+1}),R((\l x-y)^{i+1}))]}{(a,b)}}_2\\
 &&-\frac{1}{2}{\inn{[(x,y),(y(\l x-y)^i+(\l x-y)^iy,\l x(\l x-y)^i+\l(\l x-y)^ix)]}{(a,b)}}_2\\
 &&-\frac 12\inn{[\l x-y,(\l x-y)^i]}{R(B)}\\
 &&-\frac{1}{2}\inn{[x,\l(\l x-y)^i]}{yb+by}+\frac{1}{2}\inn{[y,(\l x-y)^i]}{xa+ax},
\end{eqnarray*}
Moreover according to the  $\Ad$-invariance of $P_i$ and the property  
$\inn{x}{yz}=\inn{xy}{z}, \forall$ $x$,$y$,$z\in~\g$, we obtain
\begin{eqnarray*}
 \X_{P_i\circ\phi_\l}^{Q}(x,y)[K]
&=&- {\inn{[(x,y),(R((\l x-y)^{i+1}),R((\l x-y)^{i+1}))]}{(a,b)}}_2\\
&&-\frac{1}{2}{\inn{[(x,y),(y(\l x-y)^i+(\l x-y)^iy,\l x(\l x-y)^i+\l(\l x-y)^ix)]}{(a,b)}}_2\\
&&-\frac{1}{2}\inn{y[x,\l(\l x-y)^i]+[x,\l (\l x-y)^i]y}{b}\\
&&+\frac{1}{2}\inn{x[y,(\l x-y)^i]+[y,(\l x-y)^i]x}{a}.
\end{eqnarray*}%
We then deduce  that
\begin{eqnarray*}
\X_{P_i\circ\phi_\l}^{Q}(x,y)&=&-[(x,y),(R((\l x-y)^{i+1}),R((\l x-y)^{i+1}))]\\
                         &&-\frac{1}{2}[(x,y),(y(\l x-y)^i+(\l x-y)^iy,\l x(\l x-y)^i+\l(\l x-y)^ix)]\\
                       &&+\frac{1}{2}(x[y,(\l x-y)^i]+[y,(\l x-y)^i]x,y[x,\l (\l x-y)^i]+[x,\l (\l x-y)^i]y).\\
\end{eqnarray*}%
 Let us compute separately the two last lines:
\begin{eqnarray*}
 \lefteqn{-\frac{1}{2}[(x,y),(y(\l x-y)^i+(\l x-y)^iy,\l x(\l x-y)^i+\l(\l x-y)^ix)]}\\
                      &+&\frac{1}{2}(x[y,(\l x-y)^i]+[y,(\l x-y)^i]x,y[x,\l (\l x-y)^i]+[x,\l (\l x-y)^i]y)\\
                      &=&-(x(\l x-y)^iy-y(\l x-y)^ix,\l y(\l x-y)^ix-\l x(\l x-y)^iy)\\
                      &=&([x,(\l x-y)^{i+1}],-\l [x,(\l x-y)^{i+1}])\\
                      &=&[(x,y),((\l x-y)^{i+1},-(\l x-y)^{i+1})].
\end{eqnarray*}%
Then we have Formula (\ref{vectorquad}).
\end{proof}
Let us compare the Hamiltonian vector fields for the quadratic $\RR$-Poisson bracket and the linear 
$\RR$-Poisson bracket. Formula 
(\ref{champh}) gives, in our case, the following expression for the Hamiltonian vector  field
 of $P_{i+1}\circ\phi_{\l}$ with respect to the linear $\RR$-Poisson structure:
\begin{equation}\label{champh1}
\X_{P_{i+1}\circ\phi_{\l}}=\frac{1}{2}(\l-1)[(x,y),((R-I)(\l x-y)^{i+1},(R+I)(\l x-y)^{i+1})].
\end{equation}   
Comparing Formulae
(\ref{champh1}) and  (\ref{vectorquad}), we obtain, for $\l\neq 1$, 
\begin{equation}\label{relquad}
\X_{P_i\circ\phi_{\l}}^{Q}(x,y)=\frac{2}{1-\l}\X_{P_{i+1}\circ\phi_{\l}}(x,y).
\end{equation}%
The relation (\ref{relquad}) 
 implies that for every $0\leqslant i\leqslant n-1$, we have
\begin{equation}\label{relquadline}
\renewcommand{\arraystretch}{1.9}
\left\{
\begin{array}{l}
\X_{F_{0,i}}^Q=2\X_{F_{0,i+1}},\\
\X_{F_{j,i}}^Q-\X_{F_{j-1,i}}^Q=2\X_{F_{j,i+1}},\qquad 1\leqslant j\leqslant i+1,\\
\X_{F_{i+1,i}}^Q=-2\X_{F_{i+2,i+1}}.\\
\end{array}
\right.
\end{equation}%
  Notice also that the 
Hamiltonian vector fields $\X_{F_{0,0}}^Q$ and $\X_{F_{1,0}}$   are precisely the equations of the $2$-Toda lattice on $\gl{n}(\C)$ as can be 
shown by specializing (\ref{relquadline}) to $i=0$ and  the first item of Remark \ref{remhamil}.

\begin{theorem}
The triplet $({\TP}',\FF'_{|_{{\TP}'}}, \PB_\RR^Q)$ is an integrable system.
\end{theorem}
\begin{proof}
The involutivity of the family ${\FF}'$ on $(\gg,{\PB}_{\RR})$ follows from item (1) in
Proposition  \ref{prop:quadra2}. Since  ${\FF}'$ is an  independent family  on ${\TP}'$ 
(see the proof of Proposition~\ref{linear-proof}), and 
since the cardinal of ${\FF}'$ is $\frac{n(n+3)}{2}$,
we have, according to \cite[Proposition  4.12]{Liv}, the inequality \begin{equation}\label{eq:rqnksense1}
\Rk({\TP}', {\PB}_{\RR}^Q)\leqslant 2({\rm dim}({\TP}')- {\rm Card}({\FF}'))= n^2+n-2.
\end{equation}%
 Moreover, according to  Formula (\ref{relquadline}),
the family of vector fields $\X_{{\FF}'}:=(\X_{F_{j,i}}^Q, 0\leqslant i\leqslant n-1 \textrm{ and } 0\leqslant j\leqslant i)$ and the family
of vector fields $ \X_{{\FF}}:=(\X_{F_{j,i}} , 1\leqslant i\leqslant n-1 \textrm{ and } 1\leqslant j\leqslant i+1)$ 
have the same rank at all points.
By choosing a point in $\TP$, we can deduce from the fact that $\FF$ is an integrable 
system on $\TP$ (Theorem \ref{principal-theorem}) that this rank is at least the cardinal of $\FF$
($=\frac{n(n+3)}{2}-2 $) minus the number of independent Casimir functions on $\TT$ for the linear bracket
$\PB_\RR$ ( $= n-1$).  
Since ${\FF}'$ is involutive, we have therefore the inequality
 $$\Rk{{\PB}_{\RR}^Q}\geqslant 2({\rm Rk}({\mathcal X}_{{\FF}'}) )=2({\rm Rk}({\mathcal X}_{\FF}) ) \geqslant 2( \frac{n(n+3)}{2}-2-n+1 ) =n^2+n-2.$$
 Together with (\ref{eq:rqnksense1}), this implies that the rank of the restriction
 to ${\TP}'$ of $\PB_\RR^Q $ is exactly $n^2+n-2$.
 The identity $ \card{{\FF}'}_{\arrowvert{\TP}'} =\dim{\TP}'-\frac{1}{2}\Rk({\TP}', {\PB^Q}_{\RR})$
 follows and completes the proof.
\end{proof}
\subsection{The relation between the $2$-Toda lattice and the Toda lattice }

In this section we show that the  Toda lattice  is a restriction of the $2$-Toda lattice.
 We begin by recalling the Liouville integrable system of the Toda lattice.
\subsubsection{The Toda lattice}

We give some definitions and properties  of the
   \emph{Toda lattice}, which will be useful
afterwards in this section.
\begin{definition}
\emph{(1)}
\emph{The phase space} $\TT$ of the  \emph{Toda lattice}  is the affine subspace of  $\g$ given by
\begin{equation}
\TT:=\g_{-1}\oplus\g_0+e, 
\end{equation}%
where $e=\sum_{i=1}^{\ell}e_i$.\\
\emph{(2)} The Toda lattice is the system of differential equations  on
$\TT$ given by the Lax equation
 \begin{equation}\label{todamv}
 \dot A=[A_+,A],
 \end{equation}%
 where   $A_+$ is the projection of  $A$ on  $\gp$.
\end{definition}
We consider the endomorphism  $R:=P_+-P_-$  of  $\g$, the  difference of the 
projections on  $\gp$ and  $\gm$. The Poisson $R$-bracket on $\g$ is  defined, for every   $F,G\in\FF(\g)$ and every  $x\in\g$, by
\begin{equation}\label{poissn-bracket}
{\pb{F,G}}_{R}:=\frac{1}{2}\inn{x}{[R\nabla_xF,\nabla_xG]+[\nabla_xF,R\nabla_xG]}.
\end{equation}%
\begin{theorem}\label{inte:toda-lattice}\cite[Section 4.1]{Per1990}
\emph{(1)} The  affine subspace  $\TT$ of  $\g$ is a Poisson submanifold of   $(\g,{\PB}_R)$.\\
\emph{(2)} Let  $H\in\FF(\g)$, defined for every  $x\in\g$ by
  $H(x)=\frac 12\inn{x}{x}$. The equation of the Hamiltonian  field
  $\X_H:={\Pb{H}}_R$ is the  equation of motion  (\ref{todamv}) of the Toda lattice.\\
\emph{(3)} Let  $P_1,\dots,P_{\ell}$ be a generating family of homogeneous polynomials of the algebra of $\Ad$-invariant functions on $\g$ of  degree respectively,
$m_1+1,\dots,m_{\ell}+1$. We define $\FF_0:=(P_1,\dots,P_{\ell})$. The triplet
  $(\FF_0,{\PB}_R,\TT)$ is  Liouville  integrable  system  and the equation of motion of the Toda lattice is
\begin{equation}\label{eqlaxtoda}
\dot A:={\Pb{P_1}}_R(A)=[(\nabla_AP_1)_+,A].
\end{equation}
\end{theorem}
\subsubsection{Restriction of  the $2$-Toda lattice and construction of the   Toda lattice}

The phase space $\TP$ of the $2$-Toda lattice decomposes as
\begin{equation}
\TP:=\g_{\leqslant -1}\times\g_{\geqslant 0}\oplus  \Delta(\g_{-1}\oplus
\g_0)+(e,e),
\end{equation}%
where $\Delta(\g_{-1}\oplus\g_0):=\{(x,x)\mid x\in\g_{-1}\oplus
\g_0\}$.
\begin{theorem}\label{restriction_Toda}
Let   $\gg$ be the Lie algebra $\g\times\g$ with ${\PB}_{\RR}$ the $\RR$-Poisson bracket,
   $\TP$ the phase space of the  $2$-Toda lattice and
${\TT}':= \Delta(\g_{-1}\oplus \g_0)+(e,e).$
\begin{enumerate}
\item[(1)] The submanifold   $\TT'$ is a  Poisson submanifold of
$(\gg,{\PB}_{\RR})$.
\item[(2)] Let  $(\TT,\PB_R)$ be  the  phase space of the Toda lattice, equipped
 with  the Poisson  $R$-bracket given as in (\ref{poissn-bracket}).  The map
\begin{equation}
\begin{array}{ccccc}
\varphi &:&({\TT},{\PB}_{R})&\to&({\TT}',{\PB}_{\RR})\\
        & & x             &\mapsto&(x,x)
\end{array}
\end{equation}%
is an  isomorphism of  integrable   systems. 
\end{enumerate}
\end{theorem}
\begin{proof}
{\em(1)} It is clear that ${\TT}'=\Delta(\g_{-1}\oplus\g_{0}+e)=\TP\cap\ggp$. According to  Proposition~\ref{K6},   $\TP$ is a  Poisson submanifold of  $(\gg,{\PB}_{\RR})$ while   $\ggp$ is a Poisson  submanifold  of
$(\gg,{\PB}_{\RR})$ because the
orthogonal of
$\gg$, which is  $\gg_+$,  is a Lie ideal  for the
bracket ${\LB}_{\RR}$; their intersection $\TT'$ is then a Poison submanifold of $(\gg,{\PB}_{\RR})$.\\
{\em(2)}   By using the coordinate functions on $\TT$ and on $\TT'$ we show that
  $\varphi$ is a Poisson
  isomorphism. Furthermore,  the functions of the integrable system  $\FF = (F_{k,i}, 1\leqslant i\leqslant
  \ell\textrm{ and }  0\leqslant k\leqslant m_i+1)$, restricted to  ${\TT}'$,
and pulled back on $\TT$ by  $\varphi $, give again  the functions of the family $\FF_0$ of the Toda lattice.
Specifically, for every  $i$ in
 $1,\dots, \ell$ and every  $0\leqslant k\leqslant m_i+1$, the functions $
F_{k,i}$ are all equal
(up to multiplicative constants) to 
the function  $P_i$. In fact,   for every  $x \in \g_0\oplus \g_{-1}$,
\begin{eqnarray*}
  P_i\circ \phi_\l ( \varphi (x)) &=& P_i( \phi_\l ( x,-x))\\
                              &=&  P_i (\l x-x)\\
                              &=&(\l-1)^{m_i+1} P_i (x),
\end{eqnarray*}
and  so 
  $$ \sum_{k=0}^{m_i+1} (-1)^{m_i+1-k}\l^kF_{k,i} (\varphi (x))= \sum_{k=0}^{m_i+1}(-1)^{m_i+1-k} \l^k C_{m_i+1}^k P_i (x)$$
and $F_{k,i}(\varphi(x))=C_{m_i+1}^kP_i(x)$. In conclusion, $\varphi:\TT\to \TT'$ is  Poisson isomorphic 
and $\varphi^*\FF=\FF_0$ (see item \emph{(3)} of Theorem \ref{inte:toda-lattice} for the definition of $\FF_0$). This  proves the claim.   
\end{proof}

\end{document}